\documentclass[10pt]{amsart}
\usepackage{amsfonts,amssymb,epsfig}
\usepackage[latin1]{inputenc}
\usepackage[all]{xy}
\usepackage{amsmath}
\usepackage{delarray}
\usepackage{amssymb}
\newtheorem{thm}{Theorem}[section]
\newtheorem{cor}[thm]{Corollary}
\newtheorem{lem}[thm]{Lemma}
\newtheorem{prop}[thm]{Proposition}
\theoremstyle{definition}
\newtheorem{defn}[thm]{Definition}
\newtheorem{example}[thm]{Example}
\theoremstyle{remark}

\numberwithin{equation}{section}

\renewcommand{\epsilon}{\varepsilon}

\newcommand{\dif}{{\rm d}}

\newcommand{\fg}{\mathfrak{g}}
\newcommand{\ft}{\mathfrak{t}}
\newcommand{\fj}{\mathfrak{j}}
\newcommand{\fl}{\mathfrak{l}}
\newcommand{\fs}{\mathfrak{s}}

\newcommand{\cG}{\mathcal G}
\newcommand{\cH}{\mathcal H}

\newcommand{\cO}{\mathcal O}

\newcommand{\cU}{\mathcal U}
\newcommand{\bbR}{\mathbb R}
\newcommand{\bbT}{\mathbb T}

\newcommand{\bbN}{\mathbb N}
\newcommand{\bbZ}{\mathbb Z}

\newcommand{\QED}{\hfill $\square$\vspace{2mm}}
\newcommand{\Proof}{{\bf Proof}}
\newcommand{\Remark}{{\it Remark}}

\newcommand{\Ad}{{\rm Ad}}

\begin{document}

\title
[Proper groupoids and momentum maps] {Proper groupoids and momentum
maps: linearization, affinity, and convexity \\  {\it Groupoïdes
propres et applications moment: linéarisation, caractère affine, et
convexité}}

\author{Nguyen Tien Zung}
\address{Labo. Emile Picard, UMR 5580 CNRS, UFR MIG, Université Toulouse III}
\email{tienzung@picard.ups-tlse.fr} \keywords{proper groupoid, momentum
map, linearization, affinity, convexity, symplectic, quasi-symplectic}

\subjclass{58H05, 57S15, 53D20}
\date{4th version, July 2006}

\begin{abstract} We show that proper Lie groupoids are locally linearizable.
As a consequence, the orbit space of a proper Lie groupoid is a
smooth orbispace (a Hausdorff space which locally looks like the
quotient of a vector space by a linear compact Lie group action). In
the case of proper (quasi-)symplectic groupoids, the orbit space
admits a natural integral affine structure, which makes it into an
affine orbifold with locally convex polyhedral boundary, and the
local structure near each boundary point is isomorphic to that of a
Weyl chamber of a compact Lie group. We then apply these results to
the study of momentum maps of Hamiltonian actions of proper
(quasi-)symplectic groupoids, and show that these momentum maps
preserve natural transverse affine structures with local convexity
properties. Many convexity theorems in the literature can be
recovered from this last statement and some elementary results about
affine maps.

\vspace{5mm}

\noindent {\sc Résumé}. {\it Nous montrons que les groupoides de Lie
propres sont localement linéarisables. En particulier, l'espace des
orbites d'un groupoide de Lie propre est localement isomorphe au
quotient d'un espace vectoriel par une action linéaire d'un groupe
de Lie compact. Dans le cas des groupoides (quasi-)sympletiques
propres, l'espace des orbites est une orbifold qui admet une
structure affine naturelle avec des propriétés de convexité locale.
Nous appliquons ces résultats à l'étude des actions hamiltoniennes
de groupoides (quasi-)symplectiques propres, et montrons que les
applications moment de telles actions ont une caractère affine et
des propriétés de convexité. Nous retrouvons plusieurs théorèmes de
convexité bien connus dans ce contexte.}
\end{abstract}

\maketitle

\centerline{Dedicated to Alan Weinstein on the occasion of his 60th
birthday}


\section{Introduction}

This paper consists of two parts. The first part is about the
linearization problem for Lie groupoids and (quasi-)symplectic
groupoids. The main result of this part is the \emph{local
linearization theorem} (Theorem \ref{thm:LPG}), which states that
any proper Lie groupoid with a fixed point is locally linearizable,
i.e. locally isomorphic to the action groupoid of a linear action of
a compact Lie group on a vector space. A consequence of this local
linearization theorem is the \emph{slice theorem} (Theorem
\ref{thm:LPG-Orbit}), which linearizes a proper Lie groupoid in a
neighborhood of an orbit under two additional conditions: that the
groupoid is source-locally trivial, and the orbit in question is a
manifold of finite type. This slice theorem was obtained by
Weinstein \cite{Weinstein-RPG2001} modulo Theorem \ref{thm:LPG}, and
is a generalization of the classical Koszul-Palais' slice theorem
for proper Lie group actions
\cite{Koszul-Slice1953,Palais-Slice1961} to the case of Lie
groupoids. Another immediate consequence of the local linearization
theorem is that the characteristic foliation on the base space of a
proper Lie groupoid is an orbit-like foliation in the sense of
Molino \cite{Molino-OrbitLike1994} with closed orbits, and the
corresponding orbit space (= space of orbits) is a smooth
\emph{orbispace} in the sense that it is a Hausdorff space which is
locally smoothly isomorphic to the quotient of a vector space by a
linear action of a compact Lie group. In the case of
\emph{symplectic} groupoids, Theorem \ref{thm:LPG} together with
some standard arguments imply that a slice of a proper symplectic
groupoid is locally isomorphic to a standard symplectic groupoid
$T^* G \rightrightarrows \fg^*$, where $G$ is a compact Lie group
and $\fg$ its Lie algebra (Theorem \ref{thm:LSG}), and the orbit
space is a manifold with boundary which looks locally like a Weyl
chamber (under the additional condition that the isotropy groups are
``coad-connected'', i.e. their coadjoint orbits are connected; if
this condition is not satisfied then the orbit space is an
orbifold). A similar result (Corollary \ref{cor:LQSG}) holds for
\emph{quasi-symplectic} groupoids in the sense of Xu
\cite{Xu-Momentum2003} (a.k.a. twisted presymplectic groupoids
\cite{BCWZ-TwistedDirac2003}).

The second part of this paper is about the convexity properties of
momentum maps in symplectic geometry. We will consider momentum maps in
the context of Hamiltonian spaces of quasi-symplectic groupoids
\cite{Xu-Momentum2003}, an approach which unifies the classical theory of
equivariant momentum maps for Hamiltonian group actions, Lu's momentum map
theory for actions of Poisson-Lie groups \cite{Lu-Momentum1991}, and also
Alekseev-Malkin-Meinrenken's theory of group-valued momentum maps
\cite{AMM-GroupMoment1998}. Actually, what we want to emphasize in this
paper is not the convexity, but the \emph{affinity} of momentum maps. More
precisely, we will show that if $\Gamma \rightrightarrows P$ is a proper
quasi-symplectic groupoid, then $P$ together with its characteristic
(singular) foliation admits a natural transverse integral flat affine
structure (which projects to an affine structure on the orbit space), and
any Hamiltonian $\Gamma$-space also admits a natural transverse affine
structure to a singular ``coisotropic'' foliation associated to the action
of $\Gamma$; the momentum map sends the leaves of this foliation to the
orbits of $P$, and is transversally affine, i.e. it preserves the
transverse affine structure. One then recovers various known momentum map
convexity theorems from this affine property, local convexity and some
elementary results concerning affine maps between locally convex affine
spaces.

\section{Proper groupoids}

\subsection{Linearization of proper groupoids} \hfill

Let us start by formulating the linearization problem. Consider a
Lie groupoid $\cG \rightrightarrows M$. We will always denote the
source map and the target map by $s$ and $t$ respectively. Consider
an orbit $\cO$ of $\cG$ on $M$. Then the restriction $\cG_\cO := \{p
\in G \ | \ s(p), t(p) \in \cO \}$ of $\cG$ to $\cO$ is a transitive
Lie groupoid over $\cO$, and the structure of $\cG$ induces a linear
action of $\cG_\cO$ on the normal vector bundle $N_\cO$ of $\cO$ in
$M$. (This action may be defined as follows: let $g \in \cG_\cO$ and
$x \in N_\cO$ such that $s(g)$ coincides with the projection of $x$
to $\cO$. Let $\alpha$ be a parametrized curve in $M$ such that
$\alpha(0) = s(g)$ and ${d \over d\epsilon}
\alpha(\epsilon)|_{\epsilon=0}$ projects to $x$. Let $\gamma$ be a
parametrized curve in $\Gamma$ such that $\gamma(0) = g$ and
$s(\gamma) = \alpha$. Then the projection $y$ of ${d \over
d\epsilon} t(\gamma(\epsilon))|_{\epsilon=0}$ to $N_\cO$ does not
depend on the choice of $\alpha$ and $\gamma$, and by definition
$g.x = y$). The corresponding semi-direct product $\cG_\cO \ltimes
N_\cO$ is the linear model for $\cG$ in the neighborhood of $\cO$.
The question is: do there exist a neighborhood $U$ of $\cO$ in $M$
and a neighborhood $V$ of the zero section in $N_\cO$ such that the
restriction $\cG_U = \{p \in G \ | \ s(p), t(p) \in U \}$ of the
groupoid $\cG$ to $U$ is isomorphic to $(\cG_\cO \ltimes N_\cO)_V$ ?

This linearization problem is a generalization of the problem of
linearization of Lie group actions. A special case is when $\cG = G
\ltimes M$ is the action groupoid of an action of a Lie group $G$ on a
manifold $M$ with a fixed point $m \in M$. Then the linearization of $\cG$
near $m$ is similar though somewhat weaker than the linearization of the
action of $G$ on $M$ near $m$: if the action of $G$ is linearized then the
corresponding action Lie groupoid is also linearized, and conversely if
the action groupoid $G \ltimes M$ is linearized then it means that the
action of $G$ is ``orbitally linearized'', i.e. its orbits are the same as
that of a linear action, though the action of $G$ itself may still be
nonlinear.

The classical theorems of Bochner \cite{Bochner-Linearization1945}, Koszul
\cite{Koszul-Slice1953} and Palais \cite{Palais-Slice1961} say that, under
a compactness or properness condition, smooth Lie group actions can be
linearized (near a fixed point or an orbit). On the other hand, it is easy
to construct non-proper actions (for $G = \bbR$ for example) which can't
be linearized, not even orbitally. For these reasons, in this paper we
will restrict our attention to proper groupoids.

\begin{defn}[\cite{Weinstein-RPG2001}]
A Lie groupoid $\cG \rightrightarrows M$ is called {\bf proper} if $\cG$
is Hausdorff and the map $(s,t): \cG \rightarrow M \times M$ is a proper
topological map, i.e. the preimage of a compact set is compact.
\end{defn}

\Remark. By convention, the base space (= space of objects) $M$ of a Lie
groupoid $\cG \rightrightarrows M$ is always Hausdorff, but the space of
arrows $\cG$ is a not-necessarily-Hausdorff manifold. However, \emph{all
groupoids in this paper are assumed to be Hausdorff} even when not
mentioned explicitly.

For example, the action groupoid $G \ltimes M$ of a smooth action of a Lie
group $G$ on a manifold $M$ is a proper Lie groupoid if and only if the
action of $G$ on $M$ is proper, by definition.

The above properness condition has some immediate topological
consequences, which we put together into a proposition:

\begin{prop}[\cite{MoMr-Groupoids2003,Weinstein-RPG2001}]
Let $\cG \rightrightarrows M$ be a proper Lie groupoid. Then
we have: \\
i) The isotropy group $G_m = \{ p \in \cG \ | \ s(p) = t(p) = m\}$ of any
point $m \in M$ is a compact Lie group. \\
ii) Each orbit $\cO$ of $\cG$ on $M$ is a closed submanifold of $M$. \\
iii) The orbit space $M/\cG$ together with the induced topology is a
Hausdorff space. \\
iv) If $\cH$ is a Hausdorff Lie groupoid which is Morita-equivalent to
$\cG$ then $\cH$ is also proper. \\
v) If $N$ is a submanifold of $M$ which intersects an orbit $\cO$
transversally at a point $m \in M$, and $B$ is a sufficiently small open
neighborhood of $m$ in $N$, then the restriction $\cG_B = s^{-1}(B) \cap
t^{-1}(B)$ is a proper Lie groupoid which has $m$ as a fixed point.
\end{prop}

\Proof. Points i) and v) follow directly from the definition. A
sketchy proof of point iv) can be found in Chapter 5 of
\cite{MoMr-Groupoids2003}. Point ii), which was proved in
\cite{Weinstein-RPG2001}, is a corollary of point iii). Let us give
here a proof of point iii): Let $x,y \in M$ such that their orbits
are different: $\cO(x) \cap \cO(y) = \emptyset,$ or equivalently,
$s^{-1}(y) \cap t^{-1}(x) = \emptyset$. Denote by $D^z_1 \supset
D^z_2 \supset \hdots \ni z$ a series of compact neighborhoods (i.e.
compact sets which contain open neighborhoods) of $z$ in $M$, where
$z = x$ or $y$, such that $\bigcap_{n \in \bbN} D^z_n = \{z\}$. We
have
$\bigcap_{n \in \bbN} t^{-1}(D^x_n) \cap s^{-1}(D^y_n) = t^{-1}(x) \cap
s^{-1}(y) = \emptyset.$
Since $\Gamma$ is proper, the sets $t^{-1}(D^x_n) \cap s^{-1}(D^y_n)$ are
compact. It follows that there is $n \in \bbN$ such that
$t^{-1}(D^x_n) \cap s^{-1}(D^y_n) =  \emptyset$,
or equivalently,
$\cO(D^x_n) \cap \cO(D^y_n) = \emptyset$,
where $\cO(D^x_n)$ is the union of orbits through $D^x_n$. But the orbit
space of $\cO(D^x_n)$ (resp. $\cO(D^y_n)$) is a (compact) neighborhood of
$x$ (resp., $y$) in the orbit space of $M$. Thus the orbit space of $M$ is
Hausdorff. \QED

The groupoid $\cG_B$ in point v) of the above proposition is called a {\bf
slice} of $\cG$ at $m$. This notion makes sense even when $\cG$ is not
proper. Two Lie groupoids $\Gamma_1 \rightrightarrows B_1$ and $\Gamma_2
\rightrightarrows B_2$ with fixed points $m_1 \in B_1$ and $m_2 \in B_2$
are called {\bf locally isomorphic} (near $m_1$ and $m_2$) if there are
open neighborhoods $U_1$ of $m_1$ in $B_1$ and $U_2$ of $m_2$ in $B_2$
such that $(\Gamma_1 \rightrightarrows B_1)_{U_1} := \{p \in \Gamma_1 \ |
\ s(p), t(p) \in U_1\}$ is isomorphic to $(\Gamma_2 \rightrightarrows
B_2)_{U_2}$. Recall that, similarly to the case of Lie algebroids, two
arbitrary slices of a groupoid at two points lying on a same orbit are
locally isomorphic, and the local isomorphism class may be called the
transverse groupoid structure to the orbit in question.

The main result of this paper is the following theorem, which was
conjectured by Weinstein
\cite{Weinstein-Linearization2000,Weinstein-RPG2001}:

\begin{thm}
\label{thm:LPG} Any proper Lie groupoid $\Gamma \rightrightarrows B$ with
a fixed point $m \in B$ is locally isomorphic to a linear action groupoid,
namely the action groupoid of the action of the compact isotropy group $G
= G_m$ on the tangent space $V = T_m B$.
\end{thm}

\Remark. Structural maps and manifolds of Lie groupoids are usually
assumed to be $C^\infty$-smooth, but the above theorem holds for finitely
differentiable Lie groupoids as well: if $\Gamma$ is of class $C^k$ ($k =
1, 2, \hdots, \infty$) then it can be locally linearized by an isomorphism
of class $C^k$. We suspect that the $C^\omega$ version of Theorem
\ref{thm:LPG} is also true, though we don't have a proof of it.

\Remark. In the case when the isotropy group $G$ is semisimple, Theorem
\ref{thm:LPG} (and its $C^\omega$ version) follows from the corresponding
results about linearization of Lie algebroids obtained by Monnier and the
author in \cite{Zung-Levi2002,MoZu-Levi2002}. The proof of Theorem
\ref{thm:LPG} presented in the present paper uses an averaging method and
standard Banach norm estimations, and is considerably simpler than the
Kolmogorov-Nash-Moser fast convergence method used in
\cite{Zung-Levi2002,MoZu-Levi2002}. We suspect that the results of
\cite{Zung-Levi2002,MoZu-Levi2002} might lead to a generalization of
Theorem \ref{thm:LPG} (partial linearization of non-proper Lie groupoids).
Conversely, generalizations of Theorem \ref{thm:LPG} and of
Cranic-Fernandes' theorem about integrability of Lie algebroids
\cite{CrFe-Lie2001} might lead to results about (partial) linearization of
Lie algebroids and Poisson structures.

The proof of Theorem \ref{thm:LPG} will start in the next subsection. In
the rest of this subsection, we will discuss some of its important
consequences.

An immediate consequence of Theorem \ref{thm:LPG} is that if $\cG
\rightrightarrows M$ is a proper Lie groupoid, then the characteristic
singular foliation on $M$ (by the orbits of $\cG$) is an orbit-like
foliation in the sense of Molino \cite{Molino-OrbitLike1994}. In
particular, it is a singular Riemannian foliation. Moreover, the orbit
space $M/\cG$ (together with the induced topology and smooth structure
from $M$) locally looks like the quotient of a vector space by a linear
action of a compact Lie groups. (Locally, the orbit space $M/\cG$ is the
same as the orbit space of a slice $B/\cG_B$). In analogy with the fact
that orbifolds are orbit spaces of étale proper groupoids
\cite{Haefliger-Holonomie1984,MoPr-Orbifold1997}, it would be natural to
call the orbit space (or rather the stack) of a proper Lie groupoid a
(smooth) {\bf orbispace}. In the literature there are some other similar
but maybe non-equivalent notions of orbispaces.


Another direct consequence of Theorem \ref{thm:LPG} is the following {\bf
slice theorem} for Lie groupoids, which was obtained by Weinstein (Theorem
9.1 of \cite{Weinstein-RPG2001}) under the hypothesis that Theorem
\ref{thm:LPG} is true. Recall that if $\cO$ is an orbit on $M$ of a Lie
groupoid $\cG \rightrightarrows M$, then we denote by $\cG_\cO$ the
restriction of $\cG$ to $\cO$, and by $N_\cO$ the normal bundle of $\cO$
in $M$. There is a natural linear action of $\cG_\cO$ on $N_\cO$, and we
denote by $\cG_\cO \ltimes N_\cO$ the corresponding semidirect product. A
Lie groupoid $\cG \rightrightarrows M$ is said to be {\bf source-locally
trivial} if the source map $s: \cG \rightarrow M$ makes $\cG$ into a
locally trivial fibration. An orbit (or manifold) $\cO$ is called of {\bf
finite type} if there is a proper map $f: \cO \rightarrow \bbR$ with a
finite number of critical points.

\begin{thm}[Slice theorem \cite{Weinstein-RPG2001}]
\label{thm:LPG-Orbit} Let $\cG \rightrightarrows M$ be a
source-locally trivial proper Lie groupoid, and let $\cO$ be an
orbit of finite type of $\cG$ on $M$. Then there is an invariant
neighborhood $\cU$ of $\cO$ in $M$ such that the restriction
$\cG_\cU$ of $\cG$ to $\cU$ is isomorphic to the restriction of
$\cG_\cO \ltimes N_\cO$ to a tubular neighborhood of the zero
section in $N_\cO$ (and also isomorphic to $\cG_\cO \ltimes N_\cO$
itself).
\end{thm}

We will now apply Theorem \ref{thm:LPG} to the case of symplectic
and quasi-symplectic groupoids. Recall (see, e.g.,
\cite{CoDaWe-SymplecticGroupoid1987, DZ-PoissonBook2005}) that a
symplectic groupoid is a Lie groupoid $\Gamma \rightrightarrows P$,
where $\Gamma$ is equipped with a symplectic form $\omega$ such that
the graph $\Delta = \{ (p,q, p.q) \ | \ p,q \in \Gamma, s(p) =
t(q)\}$ of the product operation of $\Gamma$ is a Lagrangian
submanifold of $\Gamma \times \Gamma \times \overline\Gamma$, where
$\overline\Gamma$ means $\Gamma$ with the opposite symplectic form
$- \omega$. If $(\Gamma,\omega) \rightrightarrows P$ is a symplectic
groupoid, then there is a unique Poisson structure $\Pi$ on $P$ such
that the source map $s: (\Gamma,\omega) \rightarrow (P,\Pi)$ is
Poisson and the target map $t: (\Gamma,\omega) \rightarrow (P,\Pi)$
is anti-Poisson; the path-connected components of the orbits of
$\Gamma$ on $P$ are the symplectic leaves of $\Pi$. For example,
consider the action groupoid $G \times \fg^* \rightrightarrows
\fg^*$ of the coadjoint action of a Lie group $G$. Identify $G
\times \fg^*$ with $T^*G$ via left translations, and equip it with
the standard symplectic form. Then it becomes a symplectic groupoid,
which we will call a {\bf standard symplectic groupoid} and denote
by $T^*G \rightrightarrows \fg^*$. The corresponding Poisson
structure on $\fg^*$ is the standard linear (Lie-) Poisson
structure.

It is easy to check that any sufficiently small slice of a (proper)
symplectic groupoid is again a (proper) symplectic groupoid: the
symplectic form of the slice is the restriction of the symplectic form of
the original symplectic groupoid to the slice. A symplectic groupoid is
called proper if it is proper as a Lie groupoid.

\begin{thm}
\label{thm:LSG} Let $(\Gamma,\omega) \rightrightarrows (P,\Pi)$ be a
proper symplectic groupoid with a fixed point $m \in P$. Then it is
locally isomorphic (as a symplectic groupoid) to the standard symplectic
groupoid $T^*G \rightrightarrows \fg^*$, where $G= G_m$ is the isotropy
group of $m$. In other words, there is an invariant neighborhood $U$ of
$m$ in $P$ and a neighborhood $V$ of $0$ in $\fg^*$ invariant under the
coadjoint action such that $((\Gamma,\omega) \rightrightarrows P)_U \cong
(T^*G \rightrightarrows \fg^*)_V$.
\end{thm}

In the above theorem, the isotropy group $G$ can be disconnected. The
proof of Theorem \ref{thm:LSG} will be given in Subsection
\ref{subsection:LSG}.

Consider now a quasi-symplectic groupoid in the sense of Xu
\cite{Xu-Momentum2003} (also known as twisted presymplectic groupoid
\cite{BCWZ-TwistedDirac2003}). This is a Lie groupoid $\Gamma
\rightrightarrows P$, equipped with a 2-form $\omega$ on $\Gamma$
and a 3-form $\Omega$ on $P$, which satisfy the following four
conditions:

i) $\dif \omega = t^* \Omega - s^*\Omega$

ii) $\dif \Omega = 0$

iii) The graph $\Delta = \{ (p,q, p.q) \ | \ p,q \in \Gamma, s(p) =
t(q)\}$ of the product operation of $\Gamma$ is isotropic with respect to
the 2-form $\omega \oplus \omega \oplus (-\omega)$ on $\Gamma \times
\Gamma \times \Gamma$.

iv) Identify $P$ with its unit section $\varepsilon(P)$ in $\Gamma$.
Due to condition iii), for each point $m \in P$, the differential
$t_*$ of the target map $t$ can be restricted to a map
\begin{equation}
t_*: \ker \omega_m \cap T_m s^{-1}(m) \rightarrow \ker \omega_m \cap T_m
P\
\end{equation}
(where $\ker \omega_m$ denotes the kernel of $\omega$ at $m$), and the
condition is that this restricted map is bijective.

The first three conditions mean that $\omega + \Omega$ is a
3-cocycle in the total de Rham complex of the groupoid $\Gamma
\rightrightarrows P$ (see \cite{Xu-Momentum2003}), and the last
condition is a weak nondegeneracy condition on $\omega$. If $\omega$
is nondegenerate and $\Omega = 0$ then one gets back to the notion
of symplectic groupoids. The base space of a quasi-symplectic
groupoid is a manifold with a twisted Dirac structure (see
\cite{BCWZ-TwistedDirac2003}). It is easy to check that a
sufficiently small slice of a (proper) quasi-symplectic groupoid is
again a (proper) quasi-symplectic groupoid (we will leave it to the
reader as an exercise).

\Remark: The convention on Lie groupoids used in this paper is different
from \cite{Xu-Momentum2003}: our source map is the target map in
\cite{Xu-Momentum2003} and vice versa.

A result of Xu (\cite{Xu-Momentum2003}, Proposition 4.8) says that if
$(\Gamma \rightrightarrows P, \omega+\Omega)$ is a quasi-symplectic
groupoid, and $\beta$ is an arbitrary 2-form on $P$, then $(\Gamma
\rightrightarrows P, \omega'+\Omega')$, where $\omega' = \omega + t^*\beta
- s^* \beta$ and $\Omega' = \Omega + \dif \beta$, is again a
quasi-symplectic groupoid, and moreover it is Morita-equivalent to
$(\Gamma \rightrightarrows P, \omega+\Omega)$ -- the notion of Morita
equivalence of quasi-Hamiltonian groupoids will be recalled in Subsection
\ref{subsection:affine_structure}. This result together with Theorem
\ref{thm:LSG} immediately leads to the following:

\begin{cor} \label{cor:LQSG}
If $(\Gamma \rightrightarrows P, \omega+\Omega)$ is a proper
quasi-symplectic groupoid with a fixed point $m$, then it is locally
isomorphic to a quasi-symplectic groupoid of the type $(T^*G
\rightrightarrows \fg^*, \omega_0+ t^*\beta - s^* \beta+\dif \beta)$,
where $(T^*G \rightrightarrows \fg^*, \omega_0)$ is the standard
symplectic groupoid of the isotropy group $G = G_m$ of $m$, and $\beta$ is
a 2-form on $\fg^*$. In particular, $(\Gamma \rightrightarrows P,
\omega+\Omega)$ is locally Morita equivalent to the standard symplectic
groupoid $T^*G \rightrightarrows \fg^*$.
\end{cor}

\Proof. Since $\Omega$ is a closed 3-form, locally it is exact, $\Omega =
\dif \beta$, and we can kill it by changing $\omega$ to $\omega_0 = \omega
+ s^*\beta - t^* \beta$. In order to apply Theorem \ref{thm:LSG}, it
remains to verify that $\omega_0$ is nondegenerate. At the fixed point
$m$, the weak nondegeneracy condition is the same as the usual
nondegeneracy condition, so $\omega_0$ is nondegenerate at $m$. The
nondegeneracy of $\omega_0$ at $m$ implies the nondegeneracy of $\omega_0$
at the other points on the isotropy group $G_m$ via the compatibility
condition iii) of the definition, so $\omega_0$ is nondegenerate at $G_m$,
and hence it is nondegenerate in a sufficiently small neighborhood of
$G_m$ in $\Gamma$. In other words, if $B$ is a sufficiently small
neighborhood of $m$ in $P$ then $(\Gamma|_B \rightrightarrows B,
\omega_0)$ will be a proper symplectic groupoid. \QED

Recall that if $G$ is a connected compact Lie group, then the orbit
space of $T^*G \rightrightarrows \fg^*$, i.e. the space of coadjoint
orbits of $G$ on $\fg^*$ can be naturally identified with the
corresponding closed Weyl chamber $\ft^*_+$. (Embed $\ft^*_+ \subset
\ft^*$ in $\fg^*$ in a natural way; then each coadjoint orbit in of
$G$ in $\fg^*$ will intersect $\ft^*_+$ at exactly one point.)
Howerver, if $G$ is disconnected, then the orbit space of $T^*G
\rightrightarrows \fg^*$ is not necessarily $\ft^*_+$, but may be a
quotient of $\ft^*_+$ by a finite group action. The reason is that,
if $G$ is disconnected, then its coadjoint action on $\fg^*$ may mix
the connected coadjoint orbits (orbits of the connected part $G^0$
of $G$) by an action of $G/G^0$. For example, consider the following
disconnected double covering $G = \bbT^2 \sqcup \theta. \bbT^2$ of
$\bbT^2$, where $\theta$ is an element such that
$\theta.g.\theta^{-1} = g^{-1}\ \forall\ g \in \bbT^2$. Then the
coadjoint action of $G^0 = \bbT^2$ on $\bbR^2 = Lie(G)^*$ is
trivial, but the coadjoint action of $\theta$ on $\bbR^2$ is given
by the map $(x,y) \mapsto (-x,-y)$. The quotient space of $\bbR^2$
by the coadjoint action of $G$ is the orbifold $\bbR^2/\bbZ_2$. In
order to avoid such orbifolds, we make the following definition:

\begin{defn}\label{def:coad-connected}
We will say that a (not necessarily connected) compact Lie group is
\emph{coad-connected} if its coadjoint orbits are connected.
\end{defn}

If a compact Lie group $G$ is coad-connected, then the orbit space
of $T^*G \rightrightarrows \fg^*$ can be naturally identified with a
Weyl chamber just like in the connected case. Of course, if $G$ is
connected then it is automatically coad-connected.

\begin{cor} \label{cor:orbitspaceQSG}
If $(\Gamma \rightrightarrows P, \omega+\Omega)$ is a proper
quasi-symplectic groupoid whose isotropy groups are coad-connected,
then the orbit space $P/\Gamma$ is a manifold with locally
polyhedral boundary: locally near each point it looks like a Weyl
chamber of a compact Lie group.
\end{cor}

If in the above corollary we drop the ``coad-connected'' condition,
then we will have to replace ``manifold'' by ``orbifold'', and
``Weyl chamber'' by ``quotient of a Weyl chamber by a finite group
action''.

\subsection{The averaging process} \hfill

Let us now proceed to the proof of Theorem \ref{thm:LPG}. It will
occupy the rest of Section 2. So from now on until the end of
Section 2, we will denote by $\Gamma \rightrightarrows  B$ a proper
Lie groupoid with a fixed point $x_0 \in B$, and by $G = G_{x_0}$
the compact isotropy group of $x_0$. A simple fact already observed
by Weinstein \cite{Weinstein-RPG2001} is that, due to the
properness, any neighborhood of $x_0$ in $B$ will contain a closed
ball-like neighborhood saturated by compact orbits of $\Gamma$. By
shrinking $B$ if necessary, we can assume that $B$ is a closed ball,
the orbits on $B$ are compact, and the source map $s: \Gamma
\rightarrow B$ is a trivial fibration.

Note that Theorem \ref{thm:LPG} is essentially equivalent to the existence
of a smooth \emph{surjective} homomorphism $\phi$ from $\Gamma$ to $G$
(after shrinking $B$ to a sufficiently small invariant neighborhood of
$x_0$), i.e. a smooth map $\phi: \Gamma \to G$ which satisfies
\begin{equation}
\label{eqn:homomorphism} \phi(p.q) = \phi(p). \phi(q) \ \ \forall \ (p,q)
\in \Gamma_{(2)} := \{ (p,q) \in \Gamma \times \Gamma, s(p) = t(q) \}\ ,
\end{equation}
and such that the restriction of $\phi$ to $G = s^{-1}(x_0) \subset
\Gamma$ is an automorphism of $G$. We may assume that this automorphism is
identity.

Indeed, if there is an isomorphism from $\Gamma \rightrightarrows B$
to an action groupoid $G \ltimes U$, then the composition of the
isomorphism map $\Gamma \to G \times U$ with the projection $G
\times U \to G$ is such a homomorphism. Conversely, we have:

\begin{lem} Assume that there is a homomorphism $\phi: \Gamma \to G$, whose restriction to $G =
s^{-1}(x_0) \subset \Gamma$ is the identity map of $G$. Then
$\Gamma$ is locally linearizable.
\end{lem}

\Proof.  Shrinking $B$ to a sufficiently small invariant
neighborhood of $z$ in $B$ if necessary, we get a diffeomorphism
\begin{equation}
(\phi, s): \Gamma \to G \times B.
\end{equation}
Denote by $\theta$ the inverse map of $(\phi, s)$. Then there is an
action of $G$ on $B$ defined by $g.x = t(\theta(g,x))$, and the map
$(\phi, s)$ will be an isomorphism from $\Gamma \rightrightarrows B$
to the action groupoid $G \ltimes B$. This action groupoid is
linearizable by the classical Bochner's theorem
\cite{Bochner-Linearization1945}, implying that the groupoid $\Gamma
\rightrightarrows B$ is linearizable. \QED

The above lemma reduces the problem of linearizing $\Gamma$ to the
problem of finding a homomorphism from $\Gamma$ to $G$ which extends
the identity map of $G$. In order to find such a homomorphism, we
will use the averaging method. The idea is to start from an
arbitrary smooth map $\phi: \Gamma \to G$ such that $ \phi|_G = {\rm
Id}.$ (Recall that $G = s^{-1}(x_0) = t^{-1}(x_0)$). Then Equality
(\ref{eqn:homomorphism}) is not satisfied in general, but it is
satisfied for $p, q \in G$. Hence it is ``nearly satisfied'' in a
small neighborhood of $G = s^{-1}(x_0)$ in $\Gamma$. In other words,
if the base $B$ is small enough, then $\phi(p.q)\phi(q)^{-1}$ is
near $\phi(p)$ for any $(p,q) \in \Gamma_{(2)}$. We will replace
$\phi(p)$ by the average value of $\phi(p.q)\phi(q)^{-1}$ for $q$
running on $t^{-1}(s(p))$ (it is to be made precise how to define
this average value). This way we obtain a new map $\widehat{\phi}:
\Gamma \to G$, which will be shown to be ``closer'' to a
homomorphism than the original map $\phi$. By iterating the process
and taking the limit, we will obtain a true homomorphism
$\phi_{\infty}$ from $\Gamma$ to $G$.

Notice that the $t$-fibers of $\Gamma \rightrightarrows B$ are compact and
diffeomorphic to $G = t^{-1}(x_0)$ by assumptions. As a consequence, there
exists a smooth Haar probability system $(\mu_x)$ on $\Gamma$, i.e. a
smooth Haar system such that for each $x \in B,$ the volume of $t^{-1}(x)$
with respect to $\mu_x$ is 1. Such a Haar probability system $(\mu_x)$ can
be constructed as follows: begin with an arbitrary Haar system $(\mu_x')$
on $\Gamma$, then define $\mu = \mu'/I$ where $I$ is the left-invariant
function $I(g) = \int_{t^{-1}(t(g))}\dif \mu'_{t(g)}$. We will fix a Haar
probability system $\mu = (\mu_x)$ on $\Gamma$.

We fix an ad-invariant metric on the Lie algebra $\fg$ of $G$ and
the induced bi-invariant metric $d$ on $G$ itself. Denote by $1_G$
the neutral element of $G$. For each number $\rho > 0$, denote by
$B_{\mathfrak{g}}(\rho)$ (resp., $B_{G}(\rho)$) the closed ball of
radius $\rho$ in $\mathfrak{g}$ (resp., $G$) centered at $0$ (resp.,
$1_G$). By resizing the metric if necessary, we will assume that the
exponential map
\begin{equation}
\exp : B_{\mathfrak{g}}(1) \to B_{G}(1)
\end{equation}
is a diffeomorphism. Denote by
\begin{equation}
\log : B_{G}(1) \to B_{\mathfrak{g}}(1)
\end{equation}
the inverse of $\exp$. Define the distance $\Delta(\phi)$ of $\phi: \Gamma
\rightarrow G$ from being a homomorphism as follows:

\begin{equation}
\Delta(\phi) = \sup_{(p,q) \in \Gamma_{(2)}}
d\left(\phi(p.q).\phi(q)^{-1}.\phi(p)^{-1}, 1_G\right)  .
\end{equation}

Let $\phi: \Gamma \to G$ be a smooth map such that $\phi|_G$ is identity.
We will assume that $\Delta(\phi) \leq 1$, so that the following map
$\widehat{\phi}: \Gamma \to G$ is clearly well-defined:

\begin{equation}
\widehat{\phi}(p) = \exp\left( \int_{q \in t^{-1}(s(p))}
\log(\phi(p.q).\phi(q)^{-1}.\phi(p)^{-1}) \dif \mu_{s(p)}\right). \phi(p)\
.
\end{equation}

Since $\mu$ is invariant under left translations, by the change of
variable $r = p.q$, we can also write $\widehat{\phi}$ as:

\begin{equation}
\widehat{\phi}(p) = \exp \left( \int_{r \in t^{-1}(t(p))}
\log(\phi(r).\phi(p^{-1}.r)^{-1}.\phi(p)^{-1}) \dif \mu_{t(p)}\right).
\phi(p)\ .
\end{equation}

Due to the commutativity of the maps $\exp$ and $\log$ with the adjoint
actions, we can also write $\widehat{\phi}$ as follows:

\begin{equation}
\widehat{\phi}(p) = \phi(p). \exp \left( \int_{q \in t^{-1}(s(p))}
\log(\phi(p)^{-1}.\phi(p.q).\phi(q)^{-1}) \dif \mu_{s(p)} \right) .
\end{equation}

It is clear that $\widehat{\phi}$ is a smooth map from $\Gamma$ to
$G$, and its restriction to $G = s^{-1}(0) \subset \Gamma$ is also
identity. The proof of the following lemma, which says that when $G$
is essentially Abelian we are done, is straightforward (we will not
need this lemma in the proof of Theorem \ref{thm:LPG}, so we will
omit its proof here):

\begin{lem} With the above notations, if $G$ is essentially commutative
(i.e. the connected component of the neutral element of $G$ is
commutative) then $\widehat{\phi}$ is a homomorphism.
\end{lem}

In general, due to the non-commutativity of $G$, $\widehat{\phi}$ is not
necessarily a homomorphism, but $\Delta(\widehat{\phi})$ (the distance of
$\widehat{\phi}$ from being a homomorphism) is of the order of
$\Delta(\phi)^2$ (Lemma \ref{lem:C0-estimate}). It means that we have the
following fast convergent iterative process: starting from an arbitrary
given smooth map $\phi: \Gamma \rightarrow G$, such that $\phi|_G = {\rm
Id}$, construct a sequence of maps $\phi_n: \Gamma \to G$ by the
recurrence formula $\phi_1= \phi$, $\phi_{n+1} = \widehat{\phi_n}$.

In the next subsections we will show that this sequence is well-defined
(after shrinking $B$ once to a smaller invariant neighborhood of $x_0$ if
necessary), and that
\begin{equation}
\phi_\infty = \lim_{n \to \infty} \phi_n
\end{equation}
exists, is smooth, and is a homomorphism from $\Gamma$ to $G$.

\Remark. The above iterative averaging process is inspired by a
similar process which was employed by Grove, Karcher and Ruh in
\cite{GrKaRu-Group1974} to prove that near-homomorphisms between
compact Lie groups can be approximated by homomorphisms. The idea of
using Grove--Karcher--Ruh's iterative averaging method was proposed
by Weinstein
\cite{Weinstein-Invariant2000,Weinstein-Linearization2000,Weinstein-RPG2001},
though he looked at the ``wrong'' map: he considered
near-homomorphisms from $G$ to the group of bisections of $\Gamma
\rightrightarrows B$ instead of near-homomorphisms from $\Gamma$ to
$G$, and was not able to prove the convergence of a corresponding
iterative averaging process.

\subsection{Spaces of maps and $C^k$-norms} \hfill

This is an auxiliary subsection where we fix some notations and write down
some standard useful inequalities.

For each $n \in \bbN$, the space of composable $n$-tuples
\begin{equation}
\Gamma_{(n)} := \{(p_1,\hdots,p_n)\ \in \Gamma \times \hdots \times \Gamma
\ |\ s(p_i) = t(p_{i+1})\ \forall\ i < n\}
\end{equation}
is smoothly diffeomorphic to $B \times G \times \hdots \times G$
($n$ copies of $G$). To fix the norms, we will fix such a
diffeomorphism (of the same smoothness class as $\Gamma$) for each
$n$. We will mainly use the manifolds $B$, $\Gamma = \Gamma_{(1)}$,
$\Gamma_{(2)}$ and $\Gamma_{(3)}$. To fix the norms on $B$ (i.e. for
maps from and to $B$), we will assume that $B$ is a closed ball
centered at $z$ in a given Euclidean space. (We will shrink $B$
whenever necessary, but the norm of the Euclidean space which
contains it will not be changed).

If $V_1$ and $V_2$ are two nonnegative numbers which depend on several
variables and parameters, then we will write $V_1 \preceq V_2$ (read $V_1$
is smaller than $V_2$ up to a multiplicative constant) if there is a
positive constant $C$ (which does not depend on the variables of $V_1$ and
$V_2$, though it may depend some some fixed parameters) such that $V_1
\leq C V_2$. We can also write $V_1 = O(V_2)$ using Landau notation. We
will write $V_1 \approx V_2$ if $V_1 \preceq V_2 \preceq V_1$.

We are interested in the $C^k$-topology ($k \leq m$ if $\Gamma$ is
only $C^m$-smooth) of the spaces of maps from $B,
\Gamma,\Gamma_{(2)},\Gamma_{(3)}$ to $\mathfrak{g}$ and $G$. We will
use $\|.\|_k$ to denote a chosen $C^k$-norm on the vector space of
$C^k$-functions from $N$ to $\fg$, where $N$ denotes one of the
spaces $B, \Gamma,\Gamma_{(2)},\Gamma_{(3)}$, etc. (Recall that $B$
is a closed ball, so $N$ is compact with boundary. It doesn't matter
much which $C^k$-norm we choose, as long as it is a norm which
provides the $C^k$-topology). We will need the following
inequalities:

If $f_r$ is a family of functions from $N$ to $\mathfrak{g}$ which
depends on a parameter $r$ which lives in a probability space $R$
then, due to the triangular inequality for a norm, we have
\begin{equation}
\label{eqn:Ck_integral} \Big\|\int_R f_r dr \Big\|_k \leq \sup_{r
\in R} \|f_r\|_k \ .
\end{equation}

 For $f_1, f_2: N \to \mathfrak{g}$, assuming that
 $\log (\exp(f_1).\exp(f_2))$ is well-defined, we have

\begin{equation}
\label{eqn:C0-exp} \|\log(\exp(f_1).\exp(f_2)) - f_1 - f_2 \|_0
\preceq \|f_1\|_0\|f_2\|_0 \ ,
\end{equation}
and
\begin{equation}
\label{eqn:C0-exp2} \|\log(\exp(f_1).\exp(f_2))\|_0 \preceq
\|f_1+f_2\|_0 \ .
\end{equation}

If, moreover, $\|f_1\|_{k-1}, \|f_2\|_{k-1} \preceq 1$ (for some
fixed $k \geq 1$) then
\begin{multline}
\label{eqn:Ck-exp} \|\log(\exp(f_1).\exp(f_2)) - f_1 - f_2 \|_k \preceq \\
\preceq \|f_1\|_0\|f_2\|_k + \|f_2\|_0\|f_1\|_k + \|f_1\|_{k-1} +
\|f_2\|_{k-1} \ ,
\end{multline}
and
\begin{multline}
\label{eqn:Ck-exp2} \|\log(\exp(f_1).\exp(f_2))\|_k \preceq \\
\preceq \|f_1+f_2\|_k + \|f_1\|_0\|f_2\|_k + \|f_2\|_0\|f_1\|_k +
\|f_1\|_{k-1} +  \|f_2\|_{k-1} \ .
\end{multline}

The above inequalities follow directly from the Leibniz rule of
derivation and the fact that the map $K: \fg \times \fg \rightarrow
\fg$ defined by $K(x,y) = \log (\exp(x).\exp(y)) - x - y$ (this map
is well-defined in a neighborhood of the origin in $\fg \times \fg$)
is an analytic map with the following properties: the linear part of
$K(x,y)$ is trivial, the quadratic part of $K$ contains only terms
of the type $x_iy_j$ (where $(x_i)$ and $(y_i)$ are the coordinates
of $x$ and $y$ respectively), and $K(x,y) = 0$ on the subspaces
$\{x=0\}$, $\{y=0\}$ and $\{x+y = 0\}$. Let us prove, for example,
Inequality (\ref{eqn:Ck-exp}). We mainly have to estimate terms of
the type $\partial^k_z K (f_1(z),f_2(z))$, where $z$ denotes an
element of $N$, and $\partial^k_z$ denotes a $k$-times partial
derivative in $z$ (in some local coordinate system). By the Leibniz
rule, $\partial^k_z K (f_1(z),f_2(z))$ is a sum of terms of the
following types: i) $\partial_{f_1} K(f_1,f_2)
\partial^{k}_z f_1$; ii) $\partial_{f_2} K(f_1,f_2) \partial^{k}_z
f_2$; iii) products of partial derivatives which are of order $\leq
k-1$ in $z$. Note that $|\partial_{x} K(x,y)| \preceq |y|$ because
$\partial_{x} K(x,0) = 0$, and therefore the terms of Type i) can be
majored by $\|f_2\|_0\|f_1\|_k$ (up to a multiplicative constant).
Similarly, the terms of Type ii) can be majored by
$\|f_1\|_0\|f_2\|_k$. On the other hand, in a product of Type iii),
we can major one factor by $\|f_1\|_{k-1} +  \|f_2\|_{k-1}$ and the
other factors by constants (because we assumed that $\|f_1\|_{k-1},
\|f_2\|_{k-1} \preceq 1$), so every term of Type iii) can be majored
by $\|f_1\|_{k-1} +  \|f_2\|_{k-1}$. Summing up all the above terms
together, we obtain Inequality (\ref{eqn:Ck-exp}).

The space $C^k(N,G)$ of $C^k$-maps from $N$ to $G$ is not a vector
space, but rather a Banach Lie group modeled on $C^k(N, \fg)$. We
will denote by $1_G$ the neutral element of this group, i.e. the map
from $N$ to $G$ which sends every element of $N$ to the neutral
element of $G$, also denoted by $1_G$. If  $f: N \rightarrow G$ is a
continuous function then we put
\begin{equation}
\|f\|_0 := \sup_{x \in N} d(f(x), 1_G),
\end{equation}
where $d(.,.)$ is the metric on $G$. If $f: N \rightarrow G$ is a
$C^k$-map ($k \geq 1$) such that $\|f\|_0 \leq 1$ (so that $\log
(f): N \rightarrow \fg$ is well defined), we put
\begin{equation}
\|f\|_k := \| \log(f) \|_k
\end{equation}
and call it the $C^k$-norm of $f$ by abuse of language. What it
measures is a $C^k$-distance from $f$ to the neutral map $1_G$. In
particular, $\|1_G\|_k = 0 \ \forall k$. We will not need $\|f\|_k$
($k \geq 1$) when $\|f\|_0 > 1$, but let us put $\|f\|_k := 2$ ($k
\geq 1$) whenever $\|f\|_0 > 1$, so that the expression $\|f\|_k$
makes sense for all $f \in C^k(N,G)$.

Let us write down some other useful standard inequalities, whose
proof is similar, if not simpler, to the proof of Inequality
(\ref{eqn:Ck-exp}).

If $\chi$ is a given smooth map from $N$ to $N'$ (e.g., the product map
from $\Gamma_{(2)}$ to $\Gamma$) and $f$ is a map from $N'$ to
$\mathfrak{g}$ or $G$, then we have (for each fixed nonnegative integer
$k$ which does not exceed the smoothness class of the groupoid):
\begin{equation}
\label{eqn:composition_est} \|f\circ \chi \|_k \preceq \|f\|_k \ .
\end{equation}

If $f_1, f_2$ are two functions from $N$ to $G$ such that
$\|f_1\|_0,\|f_2\|_0 \leq 1$ then
\begin{equation}
\label{eqn:C0_product} \| f_1.f_2 \|_0 \preceq  \|f_1\|_0 + \|f_2\|_0  \ ,
\end{equation}
and if moreover $\|f_1\|_{k-1}, \|f_2\|_{k-1} \preceq 1$ (for some fixed
$k \geq 1$), then we have:
\begin{equation}
\label{eqn:Ck_product} \| f_1.f_2 \|_k \preceq  \|f_1\|_k + \|f_2\|_k
\end{equation}
and (more refined inequalities)
\begin{equation}
\label{eqn:Ck_product2} \| f_1.f_2 \|_k - \|f_2\|_k \preceq  \|f_1\|_k +
\|f_1\|_0\|f_2\|_k  \ ,
\end{equation}
\begin{equation}
\label{eqn:Ck_conjugate} \| f_1.f_2.f_1^{-1} \|_k \preceq  \|f_2\|_k +
\|f_1\|_k. \|f_2\|_0  \ .
\end{equation}

Finally, we will need the following result about Cauchy sequences in
Banach Lie groups, which we will state as a lemma:

\begin{lem} \label{lem:CauchySequenceBanachLie}
Suppose that $f_n$ ($n\in \bbN$) are maps from $N$ to $G$ such that
$\sum_{n=1}^\infty \|f_n\|_k$ converges (for some given nonnegative
integer $k$ which does not exceed the smoothness class of $N$). Then
the product $f_n.f_{n-1}. \hdots .f_{1}$ converges in the
$C^k$-topology when $n \to \infty$ to a $C^k$-map from $N$ to $G$.
\end{lem}

\subsection{$C^0$ estimates} \hfill

\begin{lem}
\label{lem:C0-estimate} For any $\phi: \Gamma \to G$ with $\Delta(\phi)
\leq 1$ we have
\begin{equation}
\Delta(\widehat{\phi}) \preceq (\Delta(\phi))^2 \ .
\end{equation}
In particular, there is a positive constant $C_0 > 0, C_0 \leq 1$ such
that if $\Delta(\phi) \leq C_0$ then $\widehat{\phi}$ is well-defined and
\begin{equation}
\Delta(\widehat{\phi}) \leq (\Delta(\phi))^2 / C_0 \leq \Delta(\phi) \ .
\end{equation}
\end{lem}

Proof: Denote

\begin{equation}
\psi(p,q) = \phi(p.q).\phi(q)^{-1}. \phi(p)^{-1} \ ,
\end{equation}
and
\begin{equation}
\widehat{\psi}(p,q) =
\widehat{\phi}(p.q).\widehat{\phi}(q)^{-1}.\widehat{\phi}(p)^{-1} \ .
\end{equation}

Then $\psi$ and $\widehat{\psi}$ are functions from $\Gamma_{(2)}$ to $G$.
By definition of $\widehat{\phi}$, we have

\begin{equation} \label{eqn:Psihatpq}
\begin{array}{lll}
 \widehat{\psi}(p,q)  & =
  & \exp(\int_{r \in T} \log(\psi(p.q,r)) \dif\mu). \phi(p.q). \phi(q)^{-1}. \\
& & \exp(\int_{r \in T}^{-1} \log(\psi(q,r)) \dif \mu)^{-1}. \\
& & \phi(p)^{-1} \exp(\int_{r' \in t^{-1}(s(p))} \log(\psi(p,r')) \dif\mu)^{-1} \\
& = & \phi(p.q).  \phi(q)^{-1}.\phi(p)^{-1}. E(p,q) = \psi(p,q). E(p,q) \
,
\end{array}
\end{equation}
where $T = t^{-1}(s(q))$ and
\begin{equation} \label{eqn:Epq}
\begin{array}{lll}
E(p,q) & =  & Ad_{\psi(p,q)^{-1}}\exp(\int_{r \in T(s(q))}
\log(\psi(p.q,r))
d\mu). \\
& &   Ad_{\phi(p)}\exp(\int_{r \in T(s(q))}^{-1} \log(\psi(q,r))
\dif\mu)^{-1}. \\
& & \exp(\int_{r' \in t^{-1}(s(p))} \log(\psi(p,r')) \dif\mu)^{-1} \\
& = & \exp(\int_{r \in T} \log(\psi(p,q)^{-1}.\psi(p.q,r).\psi(p,q)) \dif\mu). \\
& & \exp(\int_{r \in T} \log(\phi(p).\psi(q,r)^{-1}.\phi(p)^{-1}) \dif\mu). \\
& & \exp(\int_{r \in T} \log(\psi(p,q.r)^{-1}) \dif\mu) \ \
({\rm we \ replaced} \ r' \ {\rm by} \ r = q^{-1}.r') \\
& = & \exp(\int_{r \in T} \log(A_1) \dif\mu). \exp(\int_{r \in T}
\log(A_2) \dif\mu). \exp(\int_{r \in T} \log(A_3) \dif\mu) \ ,
\end{array}
\end{equation}
where

\begin{equation}
\label{eqn:A1A2A3}
\begin{array}{l}
A_1 = \psi(p,q)^{-1}.\psi(p.q,r).\psi(p,q) \ , \\
A_2 = \phi(p).\psi(q,r)^{-1}.\phi(p)^{-1} \ , \\
A_3 = \psi(p,q.r)^{-1} \ .
\end{array}
\end{equation}

One verifies directly that

\begin{equation}
A_1.A_2.A_3 = \psi(p.q)^{-1} \ .
\end{equation}

Consider $A_1,A_2,A_3$ as maps from $\Gamma_3$ to $G$. By definition,
$\Delta(\phi) = \|\psi\|_0$. The inequality $\Delta(\phi) \leq 1$ in the
hypothesis of Lemma \ref{lem:C0-estimate}, together with the fact that the
metric on $G$ is bi-invariant, implies that

\begin{equation}
\label{eqn:AAA} \|A_1\|_0 = \|A_2\|_0 = \|A_3\|_0 = \|\psi\|_0 =
\Delta(\phi) \leq 1
\end{equation}

Applying Inequalities (\ref{eqn:C0-exp}), (\ref{eqn:Ck_integral}) and
(\ref{eqn:AAA}) several times to $E(p,q)$, we get:

\begin{equation}
\begin{array}{lll}
\log E(p,q) & =  & \epsilon_1 + \int_{r \in T} \log(A_1) \dif\mu + \int_{r
\in T} \log(A_2) \dif\mu +
 \int_{r \in T} \log(A_3) \dif\mu \\
& = & \epsilon_1+  \int_{r \in T} [\log(A_1)+ \log(A_2) + \log(A_3)] \dif\mu \\
& = & \epsilon_1 +\epsilon_2 + \int_{r \in T} \log(A_1A_2A_3) \dif\mu \\
& = & \epsilon_1 + \epsilon_2 - \log (\psi(p,q))
\end{array}
\end{equation}
where $T = t^{-1}(s(q))$ and $\epsilon_1$ and $\epsilon_2$ are some
functions such that

\begin{equation}
\|\epsilon_1\|_0, \|\epsilon_2\|_0 \preceq \Delta(\phi)^2
\end{equation}

In other words, we have $\| \log (\psi(p,q)) + \log E(p,q)\|_0 \preceq
\Delta(\phi)^2$, which implies, by Inequality (\ref{eqn:C0-exp2}), that
$\|\psi.E\|_0 \preceq \Delta(\phi)^2$. But we have $\widehat{\psi}(p,q) =
\psi(p,q). E(p,q)$, therefore

\begin{equation}
\Delta(\widehat{\phi}) = \|\widehat{\psi}\|_0 =\|\psi.E\| \preceq
\Delta(\phi)^2
\end{equation}
\hfill $\square$

Lemma \ref{lem:C0-estimate} immediately implies the uniform
convergence (i.e. convergence in the $C^0$ topology) of the sequence
of maps $\phi_n: \Gamma \to G$, defined iteratively by $\phi_{n+1} =
\widehat{\phi_n}$, beginning with an arbitrary smooth map $\phi_1$
which satisfies the inequality $\Delta(\phi_1) \leq C_0/4$. (This
inequality can always be achieved by shrinking $B$ if necessary).
Indeed, since $\|\psi_2\|_0 = \Delta(\phi_2) \leq
(\Delta(\phi_1))^2/C_0 \leq \Delta(\phi_1) \leq C_0/4 \leq 1/4$ by
Lemma \ref{lem:C0-estimate}, where
\begin{equation}
\psi_n(p,q) = \phi_n(p.q).\phi_n(q)^{-1}. \phi_n(p)^{-1} \ ,
\end{equation}
we can define $\phi_2 = \widehat{\phi_1}$, and so on, hence $\phi_n$ is
well defined for all $n \in \bbN$. By recurrence on $n$, one can show
easily  that we have
\begin{equation}
\label{eqn:C0_est} \|\psi_n\|_{0} \leq C_{0}. (b_0)^{2^{n}}  \ \ \forall n
\in \bbN, \ \ {\rm where} \ \ b_0 = \frac{1}{2} < 1 \ ,
\end{equation}
which implies in particular that $\sum_{n=1}^\infty\|\psi_n\|_0 < \infty$
(this is a very fast converging series). Put
\begin{equation}
\label{eqn:Psi_n} \Psi_n(p) = \exp \Big( \int_{q \in t^{-1}(s(p))}
\log(\psi_n(p.q)) \dif \mu \Big) \ .
\end{equation}

Then $\|\Psi_n\|_0 \preceq \|\psi_n\|_0$ (by Inequalities
(\ref{eqn:composition_est}) and (\ref{eqn:Ck_integral})), which together
with $\sum_{n=1}^\infty\|\psi_n\|_0 < \infty$ implies  that
\begin{equation}
\sum_{n=1}^\infty \|\Psi_n\|_0 < \infty \ .
\end{equation}

This last inequality implies the convergence of the product
$\Psi_n.\Psi_{n-1}\cdots\Psi_1$ in the $C^0$-topology when $n \to
\infty$. But

\begin{equation}
\Psi_n.\Psi_{n-1}\cdots\Psi_1 = \phi_{n+1}.\phi_1^{-1} \ .
\end{equation}

Thus $\phi_n$ converges in the $C^0$-topology when $n \to \infty$.
Denote by $\phi_\infty$ the limit
\begin{equation}
\phi_\infty = \lim_{n \to \infty} \phi_n \ .
\end{equation}
Then $\phi_\infty$ is a continuous homomorphism from $\Gamma$ to $G$. It
is also clear that the restriction of $\phi_\infty$ to $G$ is the identity
map from $G$ to itself.

It remains to show that $\phi_\infty$ is smooth. This is the purpose
of the next subsection, where we will show that for any $k \in
\bbN$, $k \leq m$ if $\Gamma$ belongs to the class $C^m$ only, we
have $\phi_\infty = \lim_{n \to \infty} \phi_n$ in the
$C^k$-topology as well.

\subsection{$C^k$ estimates}\hfill

Roughly speaking, we want to estimate $\psi_n$ in order to show
that, if $k$ does not exceed the smoothness class of the groupoid
$\Gamma \rightrightarrows B$, then $\sum_{n=1}^\infty \|\psi_n\|_k <
\infty$. If this series converges, then similarly to the previous
subsection, we also have $\sum_{n=1}^\infty \|\Psi_n\|_k < \infty$
where $\Psi_n = \phi_{n+1}.\phi_n^{-1}$ is given by formula
(\ref{eqn:Psi_n}), hence the product $\Psi_n.\Psi_{n-1}\hdots\Psi_1$
converges in the $C^k$-topology when $n \to \infty$, implying that
$\phi_n \to \phi_\infty$ in the $C^k$-topology.

\begin{lem}
\label{lem:Ck-estimate} Let $k \in \bbN$ be a natural number which does
not exceed the smoothness class of the groupoid $\Gamma \rightrightarrows
B$. Assume that $\|\psi\|_0 = \Delta(\phi) \leq 1$ and $\|\phi\|_{k-1}
\preceq 1$. Then we have:
\begin{equation}
\|\widehat{\psi}\|_k \preceq \|\psi\|_0\|\psi\|_k + \|\psi\|_{k-1} +
\|\psi\|_0\|\phi\|_{k-1} + \|\psi\|_0^2\|\phi\|_k
\end{equation}
\end{lem}

\Proof. Assume that $\|\phi\|_{k-1} \preceq 1$ by hypothesis of Lemma
\ref{lem:Ck-estimate}. Then by Inequality (\ref{eqn:Ck_product}) and
Inequality (\ref{eqn:composition_est}), we have $\|\psi\|_{k-1} \preceq
\|\phi\|_{k-1} \preceq 1$. Let $A_1,A_2,A_3$ be the functions defined by
Equation (\ref{eqn:A1A2A3}). We want to estimate them. For $A_3 =
\psi(p,q.r^{-1})^{-1}$, using Inequality (\ref{eqn:composition_est}), we
get:

\begin{equation}
\|A_3\|_{k-1} \preceq \|\psi\|_{k-1} \preceq 1 \ \ {\rm and} \ \
\|A_3\|_{k} \preceq \|\psi\|_{k} \ .
\end{equation}

For $A_1 =  \psi(p,q)^{-1}.\psi(p.q,r).\psi(p,q)$, using Inequality
(\ref{eqn:Ck_product}) (and Inequality (\ref{eqn:composition_est})), we
also get

\begin{equation}
\|A_1\|_{k-1} \preceq \|\psi\|_{k-1} \preceq 1 \ \ {\rm and} \ \
\|A_1\|_{k} \preceq \|\psi\|_{k} \ .
\end{equation}

The estimation of $A_2  = \phi(p).\psi(q,r)^{-1}.\phi(p)^{-1}$ is more
complicated, because it involves the function $\phi$ directly. Using
Inequality (\ref{eqn:Ck_conjugate}) we get
\begin{equation}
\|A_2\|_{k-1} \preceq \|\psi\|_{k-1} + \|\phi\|_{k-1}\|\psi\|_0 \preceq 1
\end{equation}
and
\begin{equation}
\|A_2\|_{k} \preceq \|\psi\|_{k} + \|\phi\|_{k}\|\psi\|_0 \ .
\end{equation}

Applying Inequality (\ref{eqn:Ck-exp}) and the above inequalities to
$E(p,q)$, we get that

\begin{equation}
\|\epsilon_1\|_k, \|\epsilon_2\|_k \preceq (\|\psi\|_k
+\|\phi\|_k\|\psi\|_0) \|\psi\|_0 + (\|\psi\|_{k-1} +
\|\phi\|_{k-1}\|\psi\|_0)
\end{equation}

Moreover, we have

\begin{equation}
\|E\|_i \preceq \|\psi\|_i + \|\psi\|_0\|\phi\|_i \ \ \forall
i=0,\hdots,k.
\end{equation}

Now applying Inequality (\ref{eqn:Ck-exp2}) and the last two inequalities,
we get

\begin{equation}
\begin{array}{lll}
\|\widehat{\psi}\|_k & = & \|\exp(\log(\psi)).\exp(\log(E)) \|_k \\
& \preceq & \| \log(\psi) + \log(E)\|_k + \|\psi\|_0\|E\|_k +
\|\psi\|_k\|E\|_0 +
 \|\psi\|_{k-1} + \|E\|_{k-1} \\
& = & \|\epsilon_1 + \epsilon_2\|_k + \|\psi\|_0\|E\|_k +
\|\psi\|_k\|E\|_0 +
 \|\psi\|_{k-1} + \|E\|_{k-1} \\
& \preceq & \|\psi\|_0\|\psi\|_k + \|\psi\|_{k-1} +
\|\psi\|_0\|\phi\|_{k-1} + \|\psi\|_0^2\|\phi\|_k
\end{array}
\end{equation}
\hfill $\square$

\begin{lem}
\label{lem:Ck-estimate2} With the assumptions of Lemma
\ref{lem:Ck-estimate} we have:
\begin{equation}
\|\widehat{\phi}\|_k - \|\phi\|_k \preceq \|\psi\|_k + \|\psi\|_0
\|\phi\|_k
\end{equation}
\end{lem}
\Proof. Applying Inequality (\ref{eqn:Ck_product2}) to $\widehat{\phi} =
\Psi.\phi$, we get
\begin{equation}
\|\widehat{\phi}\|_k - \|\phi\|_k \preceq \|\Psi\|_k + \|\Psi\|_0
\|\phi\|_k
\end{equation}
Now replace $\|\Psi\|_0$ by $\|\psi\|_0$ and $\|\Psi\|_k$ by $\|\psi\|_k$.
\hfill $\square$

\begin{lem}
\label{lem:Ck-estimate3} Assume that $\phi_1$ is a map from $\Gamma$ to
$G$ such that $\Delta(\phi_1) < C_0/4$, and that $\phi_{n+1} =
\widehat{\phi_n}$ for any $n \in \bbN$, as in the previous subsection. Let
$k$ be a natural number which does not exceed the smoothness class of the
groupoid $\Gamma \rightrightarrows B$. Then there is a finite positive
number $D_k > 0$ and a positive number $0 < b_k < 1$, such that for any $n
\in \bbN$ the following two inequalities hold:
\begin{equation}
\label{eqn:Ck_phi} \|\phi_n\|_k \leq D_k.(1 - 2^{-n})
\end{equation}
and
\begin{equation}
\label{eqn:Ck_psi} \|\psi_n\|_k \leq D_k.(b_k)^{2^n} \ .
\end{equation}
\end{lem}

\Proof. We will prove the above lemma by induction on $k$. When $k=0$,
Lemma \ref{lem:Ck-estimate3} is already proved in the previous section
(with $b_0 = 1/2$). Let us now assume that Inequalities (\ref{eqn:Ck_psi})
and (\ref{eqn:Ck_phi}) are true at the level $k-1$ (i.e. if we replace $k$
by $k-1$). We will show that they are true at the level $k$.

We will choose an (arbitrary) number $b_k > 0$ such that $1 > b_k > b_k^2
> b_{k-1}, b_0$. (For example, one can put $b_0 = 1/2$ and then $b_k =
(b_{k-1})^{1/3}$ by recurrence). What will be important for us is that
$b_0/b_k,b_{k-1}/b_k^2$ and $b_0/b_{k}^2$ are positive numbers which are
strictly smaller than 1.

It follows from Lemma \ref{lem:Ck-estimate} and Lemma
\ref{lem:Ck-estimate2} that there exist two positive numbers $c_1$ and
$c_2$ (which do not depend on $n$) such that we have, for any $n \in
\bbN$:
\begin{equation}
\label{eqn:psi_est2} \|\psi_{n+1}\|_k \leq c_1(\|\psi_n\|_0\|\psi_n\|_k +
\|\psi_n\|_{k-1} + \|\psi_n\|_0\|\phi_n\|_{k-1} +
\|\psi_n\|_0^2\|\phi_n\|_k)
\end{equation}
and
\begin{equation}
\|\phi_{n+1}\|_k - \|\phi_n\|_k \leq c_2(\|\psi_n\|_k + \|\psi_n\|_0
\|\phi_n\|_k) \ .
\end{equation}

We will now prove Inequalities (\ref{eqn:Ck_psi}) and (\ref{eqn:Ck_phi})
by induction on $n$. There exists a natural number $n_0$ such that for any
$n > 0$ we have

\begin{equation}
Q_1 := D_0 \left(\frac{b_0}{b_k}\right)^{2^n} +
\left(\frac{b_{k-1}}{b_k^2}\right)^{2^n} + D_0
\left(\frac{b_0}{b_k^2}\right)^{2^n} + D_0^2
\left(\frac{b_0}{b_k}\right)^{2^{n+1}} \leq \frac{1}{c_1}
\end{equation}
and
\begin{equation}
Q_2 := (b_k)^{2^n} + D_0 (b_0)^{2^n} \leq \frac{2^{-n-1}}{ c_2} \ .
\end{equation}

By choosing $D_k$ large enough, we can assume that Inequalities
(\ref{eqn:Ck_psi}) and (\ref{eqn:Ck_phi}) are satisfied for any $n \leq
n_0$. We will also assume that $D_k \geq D_{k-1}$. let us now show that if
Inequalities (\ref{eqn:Ck_psi}) and (\ref{eqn:Ck_phi}) are satisfied for
some $n \geq n_0$ then they are still satisfied when we replace $n$ by
$n+1$. (This is the last step in our induction process).

Indeed, for $\|\psi_{n+1}\|_k$, using Inequality (\ref{eqn:psi_est2}) and
the induction hypothesis, we get
\begin{equation*}
\begin{array}{l}
\|\psi_{n+1}\|_k \leq c_1\left(\|\psi_n\|_k \|\psi_n\|_0 +
\|\psi_n\|_{k-1} +
\|\psi_n\|_0 \|\phi_n\|_{k-1} + \|\psi_n\|_0^2 \|\phi_n\|_k \right)  \\
\leq c_1 \left( D_k(b_k)^{2^n}D_0(b_0)^{2^n} + D_{k-1}(b_{k-1})^{2^n} +
D_0 (b_0)^{2^n}
D_{k-1} + D_0^2 (b_0)^{2^{n+1}}D_k \right) \\
\leq D_k c_1 \left( (b_k)^{2^n}D_0(b_0)^{2^n} + (b_{k-1})^{2^n} + D_0
(b_0)^{2^n}
 + D_0^2 (b_0)^{2^{n+1}} \right) \\
= D_k c_1 Q_1 (b_k)^{2^{n+1}} \leq D_k (b_k)^{2^{n+1}}.
\end{array}
\end{equation*}

Similarly, for $\|\phi_{n+1}\|_k$ we have:
\begin{equation*}
\begin{array}{l}
\|\phi_{n+1}\|_k  \leq \|\phi_n\|_k + c_2\left(\|\psi_n\|_k +
\|\psi_n\|_0\|\phi_n\|_{k}
\right) \\
 \leq D_{k}(1 - 2^{-n}) + c_2 \left( D_k(b_k)^{2^n} + D_0(b_0)^{2^n}.D_k \right) \\
 \leq D_{k}(1 - 2^{-n}) + D_k c_2 [(b_k)^{2^n} + D_0(b_0)^{2^n}] \\
 = D_{k}(1 - 2^{-n})  + D_kc_2 Q_2
 \leq D_k (1 - 2^{-n}) + D_k 2^{-n-1} = D_k(1 - 2^{-n-1})
\end{array}
\end{equation*}
\QED

{\bf End of the proof of Theorem \ref{thm:LPG}}. Inequality
(\ref{eqn:Ck_psi}) is a sufficient condition for the
$C^k$-smoothness of $\phi_\infty$ (provided that $k$ does not exceed
the smoothness class of $\Gamma$), because it implies in particular
that $\sum_{n=1}^\infty \|\Psi_n\|_k \preceq \sum_{n=1}^\infty
\|\psi_n\|_k < \infty$, which in turns implies that the sequence of
maps $(\phi_n)$ converges in the $C^k$-topology, by Lemma
\ref{lem:CauchySequenceBanachLie}. Thus the homomorphism
$\phi_\infty: \Gamma \to G$ has the same smoothness class as
$\Gamma.$ \QED

\Remark. If we start with a near-homomorphism from $\Gamma$ to a compact
Lie group $H$ different from $G$, then our iterative averaging method
still yields a homomorphism from $\Gamma$ to $H$. So we get a
generalization of the cited Grove--Karcher--Ruh's result
\cite{GrKaRu-Group1974} about approximation of near-homomorphisms by
homomorphisms.

\subsection{Proof of Theorem \ref{thm:LSG}}\hfill
\label{subsection:LSG}

Recall that the linear part of the Poisson structure $\Pi$ at the fixed
point $m$ is isomorphic to the Lie-Poisson structure on $\fg^*$. Theorem
\ref{thm:LPG} allows us to linearize $\Gamma$ near $m$ without the
symplectic structure. The corresponding linear action of $G = G_m$ must be
(isomorphic to) the coadjoint action, so without losing generality we may
assume that $P$ is a neighborhood of $0$ in $\fg^*$, and the orbits on $P$
near $0$ are nothing but the coadjoint orbits (though the symplectic form
on each orbit may be different from the standard one). But then, as was
shown by Ginzburg and Weinstein \cite{GiWe-PoissonLie1992} using a
standard Moser's path argument, since $G$ is compact, the Poisson
structure on $P$ is actually locally isomorphic to the Lie-Poisson
structure of $\fg^*$. We can now apply the following proposition to finish
the proof of Theorem \ref{thm:LSG}:

\begin{prop}
\label{prop:g-star}
 If $G$ is a (not necessarily connected) compact Lie group and
$\fg$ is its Lie algebra, then any proper symplectic groupoid
$(\Gamma,\omega) \rightrightarrows U$ with a fixed point $0$ whose base
Poisson manifold is a neighborhood $U$ of $0$ in $\fg^*$ with the
Lie-Poisson structure and whose isotropy group at $0$ is $G$ is locally
isomorphic to $T^*G \rightrightarrows \fg^*$.
\end{prop}

{\bf Proof of Proposition \ref{prop:g-star}}. Without loss of generality,
we can assume that $\Gamma$ is source-locally trivial.

 We will first prove the above
proposition for the case when $G$ is connected. The Lie algebra $\fg$ can
be written as a direct sum $\fg = \fs \oplus \fl$, where $\fs$ is
semisimple and $\fl$ is Abelian. Denote by
$(f_1,\hdots,f_n,h_1,\hdots,h_m)$ a basis of linear functions on $\fg^*$,
where $f_1,\hdots,f_n$ correspond to $\fs$ and $h_1,\hdots,h_m$ correspond
to $\fl$. Then the vector fields $X_{s^*f_i},X_{s^*h_j}$ generate a
Hamiltonian action of $\fg$ on $(\Gamma,\omega)$. When restricted to the
isotropy group $G = G_0$ over the origin of $\fg^*$, the vector fields
$X_{s^*f_i},X_{s^*h_j}$ become left-invariant vector fields on $G$, and
the action of $\fg$ integrates to the right action of $G$ on itself by
multiplication on the right. Assume that the above Hamiltonian action of
$\fg$ integrates to a right action of $G$ on $\Gamma$. Then we are done.
Indeed, since the action is free on $G_0$, we may assume, by shrinking the
base space $U$, that the action is free on $\Gamma$. Then one can verify
directly that the map $(g,y) \mapsto \epsilon(\Ad^*_gy) \circ g$, $g \in
G$, $y \in U$, where $\epsilon : U \rightarrow \Gamma$ denotes the
identity section and $\circ g$ denotes the right action by $g$, is a
symplectic isomorphism between the restriction of  the standard symplectic
groupoid $G \times \fg^* \cong T^*G \rightrightarrows \fg^*$ to $U \subset
\fg^* $ and $\Gamma$.

In general, the action of $\fg$ on $\Gamma$ integrates to an action of the
universal covering of $G$ on $\Gamma$, which does not factor to an action
of $G$ on $\Gamma$ if the Abelian part of $G$ is nontrivial, i.e. $\fl
\neq 0$. So we may have to change the generators of this $\fg$ action, by
changing $h_1,\hdots,h_m$ to new functions $h_i'$ which are still Casimir
functions of $\fg^*$. Such a change of variables (leaving $f_i$ intact)
will be a local Poisson isomorphism of $\fg^*$. We want to choose $h_i'$
so that the Hamiltonian vector field $X_{s^*h_i'}$ are periodic, i.e. they
generate Hamiltonian $\bbT^1$-actions.

Note that for each $y \in U \subset \fg^*$, the isotropy group $G_y =
s^{-1}(y) \cap t^{-1}(y)$ of $\Gamma$ at $y$ admits a canonical injective
homomorphism to $G$ (via a a-priori non-symplectic local linearization of
$\Gamma$ using Theorem \ref{thm:LPG}). Denote by $\bbT^m_0$ the Abelian
torus of dimension $m$ in the center of $G$ (the Lie algebra of $\bbT^m_0$
is $\fl$). The coadjoint action of $\bbT^m_0$ on $\fg^*$ is trivial. It
follows that each isotropy group $G_y$ contains a torus $\bbT^m_y$ whose
image under the canonical injection to $G$ is $\bbT^m_0$. For each $q \in
\Gamma$, denote $\bbT^m_q = q. \bbT^m_{s(q)} = \bbT^m_{t(q)}.q$. Note that
if $r \in \bbT^m_q$ then $\bbT^m_r = \bbT^m_q$.

Choose a basis $\gamma_1,\hdots,\gamma_m$ of 1-dimensional sub-tori
$\bbT^m$. Translate them to each point $q \in \Gamma$ as above, we get $m$
curves $\gamma_{1,q},\hdots,\gamma_{m,q} \subset \bbT^m_q$ $\forall\ q \in
\Gamma$. Recall that, due to the fact that $\bbT^m_0$ lies in the center
of $G$, these curve are well-define and depend continuously on $q$.

Since $G = s^{-1}(0) = t^{-1}(0)$ is a Lagrangian submanifold of
$\Gamma$, the symplectic form $\omega$ of $\Gamma$ is exact (near
$G$) and we can write $\omega = \dif \alpha$. Define $m$ functions
$H_i$, $i=1,\hdots,m$ on $\Gamma$ via the following integral
formula, known as Arnold--Mineur formula for action functions of
integrable Hamiltonian systems \cite{Mineur-AA1937}:
\begin{equation}
H_i (q) = \int_{\gamma_{i,q}} \alpha .
\end{equation}

Denote by $\Gamma_{reg}$ the ``regular'' part of $\Gamma$, i.e. the set of
points $q \in \Gamma$ such that $t(q)$ is a regular point of the Poisson
structure $\Pi$ in $P$. Then $\Gamma_{reg}$ admits a natural
symplectically complete foliation by isotropic submanifolds $K_q =
q.G_{s(q)}$, and since $\Gamma$ is proper, these submanifolds are compact.
So this foliation may be viewed as the foliation by invariant tori of a
proper non-commutatively integrable Hamiltonian system. Since
$\gamma_{i,q} \subset K_q\ \forall i$, it follows from the classical
Arnold--Liouville--Mineur theorem on action-angle variables of integrable
Hamiltonian systems that $H_i$ are action functions, i.e. the Hamiltonian
vector fields $X_{H_i}$ are periodic (of period 1) and generate
$\bbT^1$-actions, and they are tangent to the isotropic submanifolds $K_q,
q \in \Gamma$. This fact is true in $\Gamma_{reg}$, which is dense in
$\Gamma$, so by continuity it's true in $\Gamma$.

By construction, the action functions $H_i$ are invariant on the leaves of
the dual coisotropic foliation of the foliation by $K_g$, $g \in
\Gamma_{reg}$, so they project to (independent) Casimir functions on $P$.
In other words, we have $m$ independent Casimir functions
$h_1',\hdots,h_m'$ such that $s^*h_i' = t^*h_i' = H_i$.

The infinitesimal action of $\fg$ on $\Gamma$ generated by Hamiltonian
vector fields $X_{s^*f_i}$, $X_{s^*h_j'}$, where the functions
$f_1,\hdots,f_n$ are as before, now integrates into an action of $S \times
\bbT^m$ on $\Gamma$, where $S$ is the connected simply-connected
semisimple Lie group with Lie algebra $\fs$. The group $S \times \bbT^m$
is a finite covering of $G$, i.e. we have an exact sequence $0 \rightarrow
\cG \rightarrow S \times \bbT^m \rightarrow G \rightarrow 0$, where $\cG$
is a finite group. Indeed, by construction, for every element $g \in \cG
\subset S \times \bbT^m$, the action $\phi(g)$ of $g$ on $\Gamma$ is
identity on the isotropy group $G$, and its differential at the neutral
element $e \in G \subset \Gamma$ is also the identity map of $T_e\Gamma$.
Since a finite power of $\phi(g)$ is the identity map on $\Gamma$, it
follows that $\phi(g)$ itself is the identity map. Hence the action of
$\cG$ on $\Gamma$ is trivial, and the action of $S \times \bbT^m$ on
$\Gamma$ factors to a Hamiltonian action of $G$ on $\Gamma$. The
proposition is proved for the case when $G$ is connected.

Consider now the case $G$ is disconnected. Denote by $G^0$ the
connected component of $G$ which contains the neutral element, and
by $\Gamma^0$ the corresponding connected component of $\Gamma$ (we
assume that the base $U$ is connected and sufficiently small). Then
$\Gamma_0$ is a proper symplectic groupoid over $U$ whose isotropy
group at $0$ is $G^0$. According to the above discussion, $\Gamma^0$
can be locally symplectically linearized, i.e. we may assume that
$\Gamma^0$ is symplectically isomorphic to $(T^*G^0
\rightrightarrows \fg^*)_U$ with the standard symplectic structure.
Consider a map $\phi: \Gamma \rightarrow G$, whose restriction to
the isotropy group $G = s^{-1}(0) \cap t^{-1}(0)$ is identity, and
whose restriction to $\Gamma^0$ is given by the projection $T^*G^0
\cong G^0 \times \fg^* \rightarrow G^0$ after the above symplectic
isomorphism from $\Gamma^0$ to $(G^0 \times \fg^* \rightrightarrows
\fg^*)_U$. We can arrange so that $\phi(p) = \phi(p^{-1})^{-1}$ for
any $p \in \Gamma$, and also $\phi(p).\phi(q) = \phi(p.q)$ for any
$p \in \Gamma^0, q \in \Gamma$. (This is possible because
$\phi|_{\Gamma_0}: \Gamma_0 \rightarrow G^0$ is a homomorphism).
Then the averaging process used in the proof of Theorem
\ref{thm:LPG} does not change the value of $\phi$ on $\Gamma$. By
repeating the proof of theorem \ref{thm:LPG}, we get a homomorphism
$\phi_\infty: \Gamma \rightarrow G$, which coincides with $\phi$ on
$\Gamma^0$.

Identifying $\Gamma$ with $G \times U$ via the isomorphism $p \mapsto
(\phi_\infty(p), s(p))$ as in the proof of Theorem \ref{thm:LPG}, and then
with $(T^*G \rightrightarrows \fg^*)_U$, we will assume that $\Gamma$, as
a Lie groupoid, is nothing but the restriction $(T^*G \rightrightarrows
\fg^*)_U$ of the standard symplectic groupoid $T^*G \rightrightarrows
\fg^*$ to $U \subset \fg^*$, and the symplectic structure $\omega$ on
$(T^*G)_U \cong G \times U$ coincides with the standard symplectic
structure $\omega_0$ on the connected component $(T^*G^0)_U \cong G^0
\times U$. For each $\theta \in G/G^0$, we will denote the corresponding
connected component of $G$ by $G^\theta$ and the corresponding connected
component of $\Gamma$ by $\Gamma^\theta$. We will use Moser's path method
to find a groupoid isomorphism of $\Gamma$ which moves $\omega$ to
$\omega_0$.

Let $f : U \rightarrow \bbR$ be a function on $U$. Then the Hamiltonian
vector fields $X^\omega_{s^*f}$ and $X^{\omega_0}_{s^*f}$ of $s^*f$ with
respect to $\omega$ and $\omega_0$ are both invariant under left
translations in $\Gamma$, and since they coincide in $\Gamma^0$ they must
coincide in $\Gamma$, because any element in $\Gamma$ can be
left-translated from an element in $\Gamma^0$. So we have a common
Hamiltonian vector field $X_{s^*f}$ for both $\omega$ and $\omega_0$.
Similarly, we have a common Hamiltonian vector field $X_{t^*f}$ for both
$\omega$ and $\omega_0$. It means that $i_X(\omega - \omega_0) = 0$ for
any $X \in T_p s^{-1}(s(p)) + T_p t^{-1}(t(p))$, which implies that
$\omega - \omega_0$ is a basic closed 2-form with respect to the
coisotropic singular foliation whose leaves are connected components of
the sets $s^{-1}(s(t^{-1}(t(p)))$, $p \in \Gamma$. In particular, for any
connected component $\Gamma^\theta$ of $\Gamma$, where $\theta \in G/G^0$,
there is a unique closed 2-form $\beta_\theta$ on $U$, which is basic with
respect to the foliation by the orbits of the coadjoint action of $G^0$ on
$U$, such that
\begin{equation}
\omega - \omega_0 = s^*\beta_\theta \ \ \ {\rm  on} \ \Gamma^\theta .
\end{equation}

The coadjoint action of $G$ on $U$ induces an action $\rho$ of $G/G^0$ on
the space of connected coadjoint orbits (orbits of $G^0$) on $U$: if $\cO$
is a connected coadjoint orbit on $U$, then $\rho(\theta) (\cO)$ is the
orbit $\Ad^*_{G^\theta} \cO$. Since $\Gamma$ is a symplectic groupoid with
respect to both $\omega$ and $\omega_0$, the closed 2-form $\omega -
\omega_0$ is also compatible with the product map in $\Gamma$. By
projecting this compatibility condition to $U$, we get the following
equality:
\begin{equation}
\beta_{\theta_1\theta_2} = \beta_{\theta_2} + \rho(\theta_2)^*
\beta_{\theta_1} \ \ \ \forall\ \theta_1,\theta_2 \in G/G_0 .
\end{equation}

Since the 2-forms $\beta_\theta$ are closed on $U$ which are basic with
respect to the foliation by connected coadjoint orbits (i.e. orbits of the
coadjoint action of $G^0$), we can write
\begin{equation} \beta_\theta = \dif \alpha_\theta,
\end{equation}
where $\alpha_\theta$ are 1-forms on $U$ which are also basic with respect
to the foliation by connected coadjoint orbits. Indeed, write
$\beta_\theta = \dif \hat{\alpha}_\theta$, then define $\alpha_\theta$ by
the averaging formula
\begin{equation}
\alpha_\theta =  \int_{G^0} (\Ad^*_g)^*\hat{\alpha}_\theta \dif \mu_{G^0},
\end{equation}
where $\mu_{G^0}$ is the Haar measure on $G^0$. Then $\alpha_\theta$ is
invariant with respect to the coadjoint action of $G^0$, and $\dif
\alpha_\theta = \beta_\theta$. One verifies easily that $\alpha_\theta$
must automatically vanish on vector fields tangent to the coadjoint
orbits, or otherwise $\beta_\theta$ would not be a basic 2-form.

Moreover, by averaging $\alpha_\theta$ with respect to the action of
$G/G^0$ via the formula
\begin{equation}
\alpha^{new}_\theta = {1 \over |G/G^0|} \sum_{\theta' \in G/G^0}
(\alpha_{\theta'\theta} - \rho(\theta)^*\alpha_\theta),
\end{equation}
we may assume that the 1-forms $\alpha_\theta$ satisfy the equation
\begin{equation}
\alpha_{\theta_1\theta_2} = \alpha_{\theta_2} + \rho(\theta_2)^*
\alpha_{\theta_1} \ \ \ \forall\ \theta_1,\theta_2 \in G/G^0 .
\end{equation}

Consider the vector field $Z$ on $\Gamma$ defined by
\begin{equation}
s^*\alpha_\theta = i_Z \omega = i_Z \omega_0 \ \  {\rm on}\ \
\Gamma_\theta .
\end{equation}
One verifies directly that the flow $\phi^t_Z$ of $Z$ preserves the
groupoid structure of $\Gamma$, and $\phi^1_Z$ moves $\omega$ to
$\omega_0$. \QED

\section{Momentum maps}

In this section, to avoid dealing with quotient spaces which are
orbifolds instead of manifolds (with boundary and corners), we will
only consider groupoids whose isotropy groups are coad-connected
(see Definition \ref{def:coad-connected}), even when we don't
mention it explicitly.

\subsection{Affine structure on base spaces} \hfill
\label{subsection:affine_structure}

In this subsection we will show that the orbit space $X = P/\Gamma$
of a proper (quasi-)symplectic groupoid $(\Gamma \rightrightarrows
P,\omega+\Omega)$ whose isotropy groups are coad-connected is not
only a manifold with locally polyhedral boundary (Corollary
\ref{cor:orbitspaceQSG}), but it also admits a natural (locally
flat) affine structure which makes it into an {\bf integral affine
manifold with locally convex polyhedral boundary} (the boundary may
be empty). It means that $X$ admits an atlas with charts modelled on
convex subsets with non-empty interior in $\bbR^k$, and the
transformations maps are integral affine, i.e. are given by elements
of the integral affine group $GL(k,\bbZ) \ltimes \bbR^k$. Moreover,
near every point $X$ is locally affine-equivalent to a Weyl chamber
of a compact Lie group (points in the interior of $X$ correspond to
tori while point on the boundary correspond to non-commutative
compact Lie groups). We will assume that $X = P/\Gamma$ is
connected.

Consider first the standard symplectic groupoid $T^*G \rightrightarrows
\fg^*$ of a compact Lie group $G$. In this case, the orbit space
$\fg^*/T^*G$ is the space of coadjoint orbits of $G$ on $\fg^*$ and can be
identified naturally with a Weyl chamber $\ft^*_+$ (here $\ft$ denotes a
Cartan subalgebra of $\fg$). The affine structure on $\fg^*/T^*G \cong
\ft^*_+$ is induced from the standard affine structure on $t^*_+$. There
is another equivalent definition of this affine structure on $\fg^*/T^*G$,
which is more intrinsic and can be generalized to arbitrary proper
(quasi-)symplectic groupoids. Let us do it immediately for a general
proper symplectic groupoid $(\Gamma \rightrightarrows P,\omega)$:

Recall from Corollary \ref{cor:orbitspaceQSG} that the orbit space
$P/\Gamma$ is a manifold with locally convex polyhedral boundary,
and the dimension $k = \dim P/\Gamma$ is also the rank (i.e. the
dimension of a Cartan torus) of each isotropic group of $\Gamma
\rightrightarrows P$. Moreover, the points on the boundary of
$P/\Gamma$ correspond to the points on $P$ whose isotropy groups are
essentially non-Abelian, while the points in the interior of
$P/\Gamma$ correspond to the points on $P$ whose isotropy groups are
essentially Abelian (i.e. the connected component of identity is a
torus of dimension $k$). For each point $z \in P$ which projects to
an interior point of $P/\Gamma$, denote by $\bbT^k_z$ the connected
component of its isotropy group $\Gamma_z$. Choose a basis
$(\gamma_1^z,\hdots\gamma_k^z)$ of $H_1(\bbT^k_z,\bbZ)$, and move it
continuously when $z$ moves (via the Gauss-Manin connection). Note
that $\bbT^k_z$ is an isotropic submanifold in $\Gamma$, and hence
the symplectic form $\omega$ is exact in a neighborhood of
$\bbT^k_z$. Denote by $\alpha$ a primitive of $\omega$, $\dif\alpha
= \omega$ in a neighborhood of $\bbT^k_z$ in $\Gamma$, and define
the following functions (in a neighborhood of $z$ in $P$):
\begin{equation} \label{eqn:MineurArnold}
F_i(z) = \int_{\gamma^z_i} \alpha, \ i=1,\hdots,k .
\end{equation}
(This is the same as the well-known Mineur-Arnold formula for action
functions of integrable Hamiltonian systems \cite{Mineur-AA1937}).
It is clear that these functions are independent. If we change
$\alpha$ to another primitive of $\omega$, then $F_i$ are changed by
additive constants, and the closed 1-forms $\dif F_i$ are not
changed.

It is easy to check that, due to the compatibility of $\omega$ with
$\Gamma$, and more precisely to the fact that for any $g \in G_z$, $T_g
G_z$ is symplectically orthogonal to $T_g\Gamma_\cO$ where $\cO$ denotes
the orbit through $z$, these local closed 1-forms $\dif F_i$ vanish on the
orbits of $\Gamma$ on $P$. If $z'$ is another point lying on the orbit
$\cO(z)$ of $z$ (not necessarily on the same connected component of the
orbit), then there is a unique natural way to transport the basis
$(\gamma_1^z,\hdots\gamma_k^z)$ of $H_1(\bbT^k_z,\bbZ)$ to a basis
$(\gamma_1^{z'},\hdots\gamma_k^{z'})$ of $H_1(\bbT^k_{z'},\bbZ)$ via the
action of $\Gamma$. Again, by the compatibility of $\omega$ with $\Gamma$,
this transportation moves closed 1-forms $(\dif F_1,\hdots,\dif F_k)$ near
$z$ to closed 1-forms $(\dif F_1,\hdots,\dif F_k)$ near $z'$ in a unique
natural way. These facts mean that the closed 1-forms $(\dif
F_1,\hdots,\dif F_k)$ can be extended in a natural single-valued way to
independent closed 1-forms on a neighborhood of $\cO(z)$ on $P$ and then
projected to independent closed 1-forms on a neighborhood of $z/\Gamma$ in
the orbit space $P/\Gamma$.

If we replace the basis $(\gamma_1^z,\hdots\gamma_k^z)$ by another
basis of $H_1(\bbT^k_z,\bbZ)$ (say by the holonomy obtained by
moving $z$ along a loop in $P$), then the vector-valued closed
1-form $(\dif F_1,\hdots,\dif F_k)$ is changed by a linear
transformation given by an element of $GL(k,\bbZ)$. It means that
while the projection of $\dif F_1,\hdots,\dif F_k$ to $P/\Gamma$ is
only locally defined on (the interior of) $P/\Gamma$ and depends on
the choice of a basis of $H_1(\bbT^k_z,\bbZ)$, it determines in a
unique way a natural integral affine structure on (the interior) of
$P/\Gamma$ ({\bf integral} means that the linear parts of the
transformation maps lie in $GL(k,\bbZ)$). To see that this affine
structure extends well to the boundary of $P/\Gamma$, we simply go
back to the local model $T^*G \rightrightarrows \fg^*$, invoking
Theorem \ref{thm:LSG}. In this local model, it is easy to check that
the affine structure defined intrinsically above in the interior of
$\fg^*/T^*G$ coincides with the affine structure obtained by
identifying $\fg^*/T^*G$ with $t^*_+$.

\Remark. if $G_1$ is a connected finite covering of a connected
compact Lie group $G$ then their corresponding groupoids $T^*G
\rightrightarrows \fg^*$ and $T^*G_1 \rightrightarrows \fg^*$ induce
the same affine structure on $t^*_+$ but maybe different integral
affine structures: the lattice of constant integral closed 1-forms
on $t^*_+$ coming from $T^*G_1 \rightrightarrows \fg^*$ is a
sublattice of the one coming from $T^*G \rightrightarrows \fg^*$.

The situation is similar in the case of proper quasi-symplectic
groupoids. A technical difference is that, since $\omega$ is not
closed in general, we have to replace Formula \ref{eqn:MineurArnold}
by another formula in order to define the analogs of $\dif F_i$: Let
$z : [0,1] \rightarrow P$ be a small path in $P$ which projects to
the interior of $P/\Gamma$. Denote by
$(\gamma_1^r,\hdots\gamma_k^r)$ a basis of $H^1(\bbT^k_{z(r)},\bbZ)$
which depends continuously on $r \in [0,1]$ (via the Gauss-Manin
connection). Then define $\alpha_i$ to be a unique closed 1-form
(defined in a neighborhood of $z(0)$) such that for any such small
path $z$ we have
\begin{equation} \label{eqn:affine_1form}
\int_{z(0)}^{z(1)} \alpha_i = \int_C \omega,
\end{equation}
where $C$ is a cylinder in $\bigcup_{r \in [0,1]} \bbT^k_{z(r)}$ whose
intersection with each $\bbT^k_{z(r)}$ is a simple closed curve
representing $\gamma^r_i$.  The compatibility of $\omega$ with $\Gamma$
implies that the above 1-form $\alpha_i$ is well-defined (i.e. does not
depend on the choice of $C$), vanishes on the orbits of $\Gamma$ and is
invariant under the action of $\Gamma$ in a natural sense, so that it can
be projected to $P/\Gamma$. The condition $\dif \omega = t^*\Omega -
s^*\Omega$ implies that the isotropy groups in $\Gamma$ are tangent to the
kernel of $\dif \omega$, which in turn guaranties that this 1-form
$\alpha_i$ is closed. This is the replacement for the 1-form $\dif F_i$ of
the symplectic case. The rest is absolutely similar to the symplectic
case.

Recall \cite{Xu-Momentum2003} that two quasi-symplectic groupoids
$(\Gamma^1 \rightrightarrows P^1, \omega^1+\Omega^1)$ and $(\Gamma^2
\rightrightarrows P^2, \omega^2+\Omega^2)$ are called {\bf Morita
equivalent} if there exists a quasi-Hamiltonian equivalence bimodule, i.e.
a manifold $M$ with the following properties: \\
i) $\Gamma^1$ acts on $M$ from the left with momentum map $\mu^1$,
$\Gamma^2$ acts on $M$ from the right with momentum map $\mu^2$, and the
two actions commute. Moreover, the actions of $\Gamma^1$ and $\Gamma^2$ on
$M$ are free, the momentum maps are submersions, and the orbits of
$\Gamma_1$
on $M$ are precisely the fibers of $\mu_2$ and vice versa. \\
ii) There is a 2-form $\sigma$ on $M$ which makes it into a
quasi-Hamiltonian $\Gamma^1 \times \overline{\Gamma^2}$-space, where
$\overline{\Gamma^2}$ means $(\Gamma^2 \rightrightarrows
P^2,-\omega-\Omega)$, and the (left) action of $\Gamma^1 \times
\overline{\Gamma^2}$ on $M$ is given by $(g_1,g_2). m := g_1.m. g_2^{-1}$.

When two Lie groupoids are Morita equivalent, they have the same orbit
spaces up to isomorphisms. In the case of proper quasi-symplectic
groupoids, one can check easily that their orbit spaces also have the same
integral affine structure, because we can ``move'' Formula
\ref{eqn:affine_1form} from one groupoid to another via a
quasi-Hamiltonian equivalence bimodule. Summarizing, we have:

\begin{thm}
\label{thm:affine-base} If $(\Gamma \rightrightarrows P, \omega
+\Omega)$ is a proper quasi-symplectic groupoid with coad-connected
isotropy groups, then its orbit space $P/\Gamma$ admits a natural
structure of an integral affine manifold which near each point is
locally affine-isomorphic to a Weyl chamber of a compact Lie group,
and which depends only on the Morita equivalence class of $(\Gamma
\rightrightarrows P, \omega +\Omega)$.
\end{thm}

\Remark. The affine structure on $P/\Gamma$ can be lifted to $P$ to
become a transverse affine structure to the orbits of $\Gamma$ in
$P$. It can also be lifted to $\Gamma$ to become a transverse affine
structure to the foliation in $\Gamma$ given by submanifolds
$s^{-1}(t^{-1}(m))$, $m \in P$. Of course, one has a similar
intrinsic definition for these transverse affine structures in $P$
and $\Gamma$.

\begin{example} Consider the AMM (Alekseev--Malkin--Meinrenken)
groupoid \cite{Xu-Momentum2003}: it is the action groupoid $G \times G
\rightrightarrows G$ of the conjugation action of a compact Lie group $G$,
equipped with a natural quasi-symplectic structure arising from the theory
of group-valued momentum maps \cite{AMM-GroupMoment1998}. Xu
\cite{Xu-Momentum2003} showed a natural equivalence between
quasi-Hamiltonian spaces with $G$-valued momentum maps and
quasi-Hamiltonian spaces of the AMM groupoid. The orbit space of the AMM
groupoid is naturally affine-equivalent to a Weyl alcove of $G$. In
particular it is a convex affine polytope.
\end{example}

\subsection{Affinity  and local convexity of momentum maps} \hfill

Consider a Hamiltonian action of a symplectic groupoid $(\Gamma
\rightrightarrows P, \omega)$ on  a symplectic manifold
$(M,\sigma)$, i.e. an action of $\Gamma$ on $M$ which is compatible
with the symplectic forms in the following sense (see
\cite{MiWe-Groupoid1988}): the graph $\{(g,x,g.x) \ | \ g \in
\Gamma, x \in P, s(g) = \mu(x) \}$ of the action is an isotropic
submanifold in $(\Gamma,\omega) \times (M,\sigma) \times (M,
-\sigma)$. Here $\mu$ denotes the momentum map of the action; it is
a Poisson map from $M$ to $P$. This is a generalization of
Hamiltonian actions of Lie groups, because, as was shown by Mikami
and Weinstein \cite{MiWe-Groupoid1988}, there is a natural 1-1
correspondence between Hamiltonian actions of a given Lie group $G$
with equivariant momentum maps and Hamiltonian actions of the
symplectic groupoid $T^*G \rightrightarrows \fg^*$, in the following
sense: If $T^*G \rightrightarrows \fg^*$ acts on a symplectic
manifold $(M,\sigma)$ with momentum map $\mu$, then $\mu$ is also
the equivariant momentum map of a Hamiltonian action of $G$ on $M$
defined as follows:
\begin{equation} \label{eqn:HamActions}
g. x  = (L_g \mu(x)).x ,
\end{equation}
where $g \in G$,  $x \in M$, $L_g$ means left translation by $g$ in
$T^*G$, $g.x$ means the action of $g \in G$ on $x$ and $(L_g
\mu(x)).x$ means the action of $L_g \mu(x) \in T^*G$ on $x$ (note
that $s(L_g \mu(x)) = \mu(x)$). Conversely, if $\mu$ is the momentum
map of a Hamiltonian action of $G$ on $(M,\sigma)$, then Formula
(\ref{eqn:HamActions}) defines a Hamiltonian action of $T^*G
\rightrightarrows \fg^*$ on $(M,\sigma)$ with the same momentum map.
In particular, the orbits of the action of $G$ on $(M,\sigma)$ are
the same as the orbits of the action of $T^*G \rightrightarrows
\fg^*$. Remark that $G$ can be disconnected, in which case, by a
Hamiltonian $G$-action we mean a symplectic action of $G$ on a
symplectic manifold $(M,\sigma)$ together with a $G$-equivariant
Poisson map $\mu: M \rightarrow \fg^*$ (the equivariant momentum
map).

Theorem \ref{thm:LSG}, together with the above equivalence between
Hamiltonian $G$-actions and Hamiltonian $(T^*G \rightrightarrows
\fg^*)$-actions, leads immediately to the following proposition:

\begin{prop} Denote by $\mu$ the momentum map of a Hamiltonian action of a
proper symplectic groupoid $(\Gamma \rightrightarrows P, \omega)$ on
a symplectic manifold $(M,\sigma)$. Let $m$ be an arbitrary point of
$M$, and denote by $N$ a submanifold in $P$ which intersects the
symplectic leaf $\cO(\mu(m))$ of $\mu(m)$ in $P$ transversally at
$\mu(m)$. Then there is a small neighborhood $B$ of $\mu(m)$ in $N$
with the following properties:
\\
i) $M_B := \mu^{-1}(B)$ is a symplectic submanifold of $M$ which
intersects the orbits of the action of $\Gamma$ on $M$ transversally. \\
ii) $(\Gamma_B \rightrightarrows B,\omega)$ is isomorphic to $(T^*G
\rightrightarrows \fg^*)_U$, where $G = G_{\mu(m)}$ is the isotropy group
of $\mu(m)$ and $U$ is a neighborhood of $0$ in $\fg^*$. Denote by $\phi:
B \rightarrow \fg^*$ a corresponding isomorphism from $B$ to $U \subset
\fg^*$. \\
iii) The induced Hamiltonian action of $(\Gamma_B \rightrightarrows
B,\omega)$ on $(M_B,\sigma)$ is equivalent to a Hamiltonian action of $G$
on $(M_B,\sigma)$ with the equivariant momentum map $\phi \circ \mu: M_B
\rightarrow \fg^*$.
\end{prop}

The proof is straightforward. \QED

More generally, we may consider a quasi-Hamiltonian space $(M,\sigma)$, in
the sense of Xu \cite{Xu-Momentum2003}, of a proper quasi-symplectic
groupoid $(\Gamma \rightrightarrows P,\omega+\Omega)$ . It means that
$\Gamma$ acts on $M$, and the following compatibility and weak
nondegeneracy conditions are satisfied: \\
i) $\dif \sigma = \mu^*\Omega$, where $\mu$ denotes the momentum map. \\
ii) The graph of the action is isotropic with respect to the 2-form
$\omega \oplus \sigma \oplus (-\sigma)$. \\
iii) $\forall m \in M, \ \ker \sigma_m = a_* (T_{\mu(m)}s^{-1}(\mu(m))
\cap \ker \omega_{\mu(m)})$. Here $a$ denotes the action map
$s^{-1}(\mu(m)) \rightarrow M, a(g) := g.m$. In particular, if $\omega$ is
nondegenerate then $\sigma$ is also nondegenerate.

Similarly to the case of Hamiltonian actions of proper symplectic
groupoids, Corollary \ref{cor:LQSG} leads to the following proposition,
whose proof is straightforward:

\begin{prop} Denote by $\mu$ the momentum map of a quasi-Hamiltonian space
$(M,\sigma)$ of a proper quasi-symplectic groupoid $(\Gamma
\rightrightarrows P, \omega+\Omega)$. Let $m$ be an arbitrary point
of
 $M$, and denote by $N$ a submanifold in $P$
which intersects the symplectic leaf $\cO(\mu(m))$ of $\mu(m)$ in
$P$ transversally at $\mu(m)$. Then there is a neighborhood $B$ of
$\mu(m)$ in $N$ and a primitive 2-form $\beta$ of the pull-back of
$\Omega$ to $B$ ($\dif \beta = \Omega_B)$, with the following
properties:
\\
i) $(M_B := \mu^{-1}(B),\sigma - \mu^*\beta)$ is a symplectic submanifold
of $M$ which intersects the orbits of the action of $\Gamma$ on $M$ transversally. \\
ii) $(\Gamma_B \rightrightarrows B,\omega + s^*\beta - t^*\beta)$ is a
proper symplectic groupoid which is isomorphic to $(T^*G \rightrightarrows
\fg^*)_U$, where $G = G_{\mu(m)}$ is the isotropy group of $\mu(m)$ and
$U$ is a neighborhood of $0$ in $\fg^*$. Denote by $\phi: B \rightarrow
\fg^*$ a corresponding isomorphism from $B$ to $U \subset
\fg^*$. \\
iii) The induced action of $(\Gamma_B \rightrightarrows B,\omega +
s^*\beta - t^*\beta)$ on $(M_B,\sigma  - \mu^*\beta)$ is Hamiltonian, and
is equivalent to the Hamiltonian action of $G$ on $(M_B,\sigma -
\mu^*\beta)$ associated to the momentum map $\phi \circ \mu: M_B
\rightarrow \fg^*$.
\end{prop}

The above proposition means that locally, near a level set of the momentum
map and after going to a slice, a quasi-Hamiltonian space of a proper
quasi-symplectic groupoid is the same as a Hamiltonian space of a compact
Lie group. So it is natural to expect that many results concerning
momentum maps of Hamiltonian actions of compact Lie groups apply to
quasi-Hamiltonian spaces of proper quasi-symplectic groupoids as well. We
will be interested in their local convexity properties, so let us recall
the following local convexity result in the ``classical'' setting:

Consider a Hamiltonian action of a coad-connected compact Lie group
$G$ on a symplectic manifold $(M,\sigma)$ with an equivariant
momentum map $\mu: M \rightarrow \fg^*$. We can factorize $\mu$ by
the action of $G$ to get a kind of reduced momentum map:
\begin{equation} \mu_{/G} : M/G \rightarrow \fg^*/G \cong \ft^*_+ .
\end{equation}
We will assume that $\mu_{/G}^{-1}(0) = \mu^{-1}(0) / G $ is not empty,
and denote by $N/G$ a connected component of $\mu_{/G}^{-1}(0)$. (The
subset $N$ of $M$ is not necessarily connected, but $N/G$ is connected; we
assume that $M$ itself is without boundary).

\begin{prop}[Kirwan \cite{Kirwan-Convexity1984}]
\label{prop:LocalConvexity} With the above notations, there is a
neighborhood $U/G$ of $N/G$ in $M/G$ such that $\mu_{/G} (U/G)$ is a
neighborhood of $0$ in a closed convex polyhedral cone $C$ of vertex $0$
in $\ft^*_+$, and that for any $c \in \mu_{/G} (U/G)$ the reduced level
set $\mu_{/G}^{-1}(c) \cap U/G$ is connected.
\end{prop}

The above proposition was proved by Kirwan
\cite{Kirwan-Convexity1984} using Morse theory. See also, e.g.,
\cite{HiNePl-Convexity1994,Sjamaar-Convexity1998,Knop-Convexity2002}
for additional information about local properties of momentum maps.
Strictly speaking, these papers consider only the case when $G$ is
connected, but the case when $G$ is disconnected but coad-connected
is the same because of $G$-equivariance. \QED

Consider now a quasi-Hamiltonian space $(M,\sigma)$ of a proper
quasi-symplectic groupoid $(\Gamma \rightrightarrows P,\omega+\Omega)$
with momentum map $\mu$. Factorize $\mu$ by the action of $\Gamma$ to get
the reduced momentum map
\begin{equation} \mu_{/\Gamma} : M/\Gamma \rightarrow P/\Gamma .
\end{equation}

The orbit space $M/\Gamma$ of the action of $\Gamma$ on $M$ is foliated by
the connected components of the preimages of the reduced map
$\mu_{/\Gamma}$. Denote by $\widehat{M}$ the quotient space of this
(singular) foliation (together with the induced topology). Then the
reduced momentum map $\mu_{/\Gamma}$ projects to a map
\begin{equation} \hat\mu : \widehat{M} \rightarrow P/\Gamma.
\end{equation}
We will call $\hat\mu$ the {\bf transverse momentum map}. We will show
that $\widehat{M}$ and $\hat\mu$ enjoy very good affine properties:

\begin{thm}\label{thm:affine-momentum}
Let $(M,\sigma)$ be a quasi-Hamiltonian space of a proper
quasi-symplectic groupoid $(\Gamma \rightrightarrows
P,\omega+\Omega)$ with coad-connected isotropy groups, with momentum
map $\mu$. Then with the
above notations we have: \\
i) $\widehat{M}$ is a manifold with locally convex polyhedral boundary
(the boundary can be empty). \\
ii) There is a natural affine structure of $\widehat{M}$ which makes it
into an affine manifold which near boundary points is locally
affine-isomorphic to convex polyhedral cones. \\
iii) The transverse momentum map $\hat\mu$ is locally injective and
affine (i.e. the pull-back of an affine function is affine).
\end{thm}

\Proof. The affine structure on $\widehat{M}$ can be defined either by the
pull-back via $\hat\mu$ of the affine structure on $P/\Gamma$ once we
establish that $\hat\mu$ is locally injective with locally polyhedral
image, or by local 1-forms defined by a formula similar to Formula
\ref{eqn:affine_1form}: replace $\omega$ by $\sigma$, and 1-cycles on
$\bbT^k_z$ ($z \in P)$ by 1-cycles on $\bbT^k_m = \bbT^k_z . m$ ($m \in M$
with $\mu(m) = z$). The two definitions are equivalent, and the obtained
affine structure is actually an integral affine structure. Since the
statements of the above theorem are local, we can work locally, and
replace $(\Gamma \rightrightarrows P,\omega+\Omega)$ by a standard
symplectic groupoid $T^*G \rightrightarrows \fg^*$, in view of Theorem
\ref{thm:LSG} and Corollary \ref{cor:LQSG} and a result of Xu
\cite{Xu-Momentum2003} which says that Morita equivalent quasi-symplectic
groupoids have equivalent quasi-Hamiltonian spaces, i.e. we are reduced to
the case of a Hamiltonian action of a compact Lie group on a symplectic
manifold. But then it becomes nothing more than Proposition
\ref{prop:LocalConvexity}. Details are left to the reader. \QED

\Remark. One can describe the boundary points of $P/\Gamma$ and of
$\widehat{M}$ via the degeneration of 1-cycles on tori. For example, when
a point $z$ in the interior of $P/\Gamma$ goes to a boundary point, then
the torus $\bbT^k_z$ becomes a noncommutative compact group, and (at
least) one of the 1-cycles on $\bbT^k_z$ vanishes homotopically, and the
local affine function corresponding to that cycle admits a minimum or
maximum value on the corresponding boundary face of $P/\Gamma$. For a
point $m$ of $M$ which goes to a boundary point of $\widehat{M}$ after the
projection $M \rightarrow \widehat{M}$ (we assume that $M$ is without
boundary), the situation is similar: the torus $\bbT^k_m$ collapses to a
smaller torus, either because the torus $\bbT^k_z$ (with $s(z) = \mu(m)$)
degenerates to a noncommutative group (i.e. we get a boundary point of
$\widehat{M}$ which maps to a boundary point of $P/\Gamma$), or the action
of $\bbT^k_z$ degenerates (i.e. there is a subtorus whose action becomes
trivial), or both.

\Remark. In the case of Hamiltonian torus actions, the fact that
$\widehat{M}$ is an affine manifold with locally convex polyhedral
boundary was probably first pointed out in \cite{CoDaMo-Moment1988}.

\subsection{Global convexity of momentum maps} \hfill

Theorem \ref{thm:affine-momentum} allows us to reduce the problem of
convexity of momentum maps to a problem concerning affine maps between
locally convex affine manifolds . For example, when the orbit space
$P/\Gamma$ can be affinely embedded or at least immersed into $\bbR^k$ (it
happens if $P/\Gamma$ is simply-connected), we can use the following two
simple lemmas about affine maps to prove global convexity:

\begin{lem} \label{lem:convex_affine_map1}
Let $X$ be a connected compact locally convex affine manifold (with
boundary), and $\phi: X \rightarrow \bbR^k$ a locally injective affine map
from $X$ to $\bbR^k$ with the standard affine structure. Then $\phi$ is
injective, and its image $\phi(X)$ is convex in $\bbR^k$.
\end{lem}

\Proof. The proof is elementary. Take an arbitrary point $x \in X$. For
each nontrivial vector $v \in T_x X$ (if $x$ lies on the boundary of $X$
then $T_x X$ means the convex tangent cone of $X$ at $x$), denote by $l_v$
the maximal affine segment lying in $X$ which begins at $x$ and going in
the direction of $v$. Because $\phi$ is affine locally injective, the map
$\phi|_{l_v} : l_v \rightarrow \phi(l_v)$ is affine injective, and
$\phi(l_v)$ is an affine segment in $\bbR^k$. Because $X$ is compact,
$\phi(X)$ is compact and therefore $\phi(l_v)$ and $l_v$ must be compact
too, i.e. $l_v$ is a closed bounded segment. Denote by $y_v$ the other end
of $l_v$. Then $y_v$ lies on the boundary $\partial X$ of $X$. It follows
easily from the local convexity of the boundary of $X$ that the union $
\bigcup_{v \in T_xX} l_v$ is an open subset of $X$. Note that $\bigcup_{v
\in T_xX} l_v$ is also closed in $X$. Indeed, if $y_n \in l_{v_n}$, $\lim
y_n = y \neq x$, then we can assume (after taking a subsequence of $(y_n)$
and resizing $(v_n)$) that $0 \neq \lim v_n = v \in T_x M$, and it follows
easily from the compactness that $y \in l_v$. Since $X$ is closed, we have
$X = \bigcup_{v \in T_xX} l_v$, i.e. $X$ is star-shaped with respect to
$x$. It follows from the local injectivity of $\phi$ at $x$ that $\phi$ is
in fact injective, and $\phi(X)$ is star-shaped with respect to $\phi(x)$.
But since $x$ is arbitrary, it means that $\phi(X)$ is convex. \QED

\begin{lem} \label{lem:convex_affine_map2}
Let $X$ be a connected locally convex affine manifold, and $\phi: X
\rightarrow \bbR^k$ a proper locally injective affine map from $X$ to
$\bbR^k$ with the standard affine structure. Then $\phi$ is injective, and
its image $\phi(X)$ is convex in $\bbR^k$.
\end{lem}

\Proof. Absolutely similar to the proof of Lemma
\ref{lem:convex_affine_map1}, via the fact that $X = \bigcup_{v \in T_xX}
l_v$. Another proof goes as follows. Denote by $B_n$ the closed ball of
radius $n$ centered at $\phi(x)$ in $\bbR^k$, and by $X_n$ the connected
component of $\phi^{-1}(B_n)$ which contains $x$. Then one checks easily
that $X_n \subset X_{n+1}$, $\cup_{n=1}^\infty X_n = X$, and each $X_n$ is
compact locally convex. Hence, by Lemma \ref{lem:convex_affine_map1}, the
restriction of $\phi$ to $X_n$ is injective and $\phi(X_n)$ is convex.
Since $\phi(X_n) \subset \phi(X_{n+1})$, it follows that $\phi(X) =
\cup_{n=1}^\infty \phi(X_n)$ is convex. \QED

\Remark. The above lemmas are along the lines of the ``local-global
principle'' of convexity \cite{CoDaMo-Moment1988,HiNePl-Convexity1994};
here we make this principle simpler by formulating it terms of pure affine
geometry.

Many known convexity theorems of momentum maps concern the case when
$P/\Gamma$ is contractible, and can be recovered from the above two lemmas
and Theorem \ref{thm:affine-momentum}. These include, for example:

\begin{itemize}
\item Kirwan's convexity theorem \cite{Kirwan-Convexity1984}, which was
conjectured and partially proved by Guillemin and Sternberg
\cite{GuSt-Convexity1982}: If a compact Lie group $G$ acts in a
Hamiltonian way on a connected compact symplectic manifold $M$ with
an equivariant momentum map $\mu: M \rightarrow \fg^*$, and $t^*_+$
denotes a Weyl chamber in $\fg^*$, then $\mu(M) \cap \ft^*_+$ is a
convex polytope. In this case $P/\Gamma \cong \ft^*_+$, and Kirwan's
theorem follows from Theorem \ref{thm:affine-momentum} and Lemma
\ref{lem:convex_affine_map1}. The generalization of Kirwan's theorem
to the case of non-compact symplectic manifolds with proper momentum
maps \cite{HiNePl-Convexity1994} follows from Theorem
\ref{thm:affine-momentum} and Lemma \ref{lem:convex_affine_map2}.

\item Flaschka--Ratiu's convexity theorem for momentum maps of Poisson
actions of compact Poisson-Lie groups \cite{FlRa-Convexity1996}: again
$P/\Gamma \cong t^*_+$.

\item Alekseev--Malkin--Meinrenken's convexity theorem for group-valued
momentum maps \cite{AMM-GroupMoment1998}: $P/\Gamma$ is a Weyl alcove.

\item Weinstein's convexity theorem for certain Hamiltonian actions of a
\emph{noncompact} semisimple Lie group $G$ which admits a compact Cartan
subgroup (\cite{Weinstein-Noncompact2001}, Theorem 3.3): here $P/\Gamma$
may be identified with a Weyl chamber of a compact Cartan subalgebra. A
special case of this theorem with $G = Sp(2k,\bbR)$ is a beautiful result
concerning frequencies of positive-definite quadratic Hamiltonian
functions (\cite{Weinstein-Noncompact2001}, Theorem 4.1) reproduced below
(Theorem \ref{thm:Weinstein_convexity}). The paper
\cite{Weinstein-Noncompact2001} is actually one of the original
motivations for the study of proper groupoids and momentum maps suggested
to us by Weinstein.
\end{itemize}

\begin{thm}[Weinstein \cite{Weinstein-Noncompact2001}]
\label{thm:Weinstein_convexity} For any positive-definite quadratic
Hamiltonian function $H$ on the standard symplectic space $\bbR^{2k}$,
denote by $\phi(H)$ the $k$-tuple $\lambda_1 \leq \hdots \leq \lambda_k$
of frequencies of $H$ ordered non-decreasingly, i.e. $H$ can be written as
$H = \sum \lambda_i(x_i^2 + y_i^2)/2$ in a canonical coordinate system.
Then for any two given positive nondecreasing $n$-tuples $\lambda =
(\lambda_1,\hdots,\lambda_k)$ and $\gamma = (\gamma_1,\hdots,\gamma_k)$,
the set
\begin{equation}
\Phi_{\lambda,\gamma} = \{\phi(H_1+H_2) \ | \ \phi(H_1) = \lambda,
\phi(H_2) = \gamma \}
\end{equation}
is a closed, convex, locally polyhedral subset of $\bbR^k$.
\end{thm}

\Remark. The above set $\Phi_{\lambda,\gamma}$ is closed but not bounded.
For example, when $k =1$ then $\Phi_{\lambda,\gamma}$ is a half-line.

There are cases when $P/\Gamma$ can't be immersed into $\bbR^k$. An
example is the situation of locally Hamiltonian torus actions studied by
Giacobbe \cite{Giacobbe-Convexity2000}, where $P/\Gamma$ is isomorphic the
quotient of $\bbR^k$ by a lattice. Then the image $\mu(M)/\Gamma$ of a
proper momentum map $\mu$ in $P/\Gamma$ may or may not be locally convex,
due to overlapping. For example, one can easily construct a convex
polytope in $\bbR^k$ such that its image under the projection from
$\bbR^k$ to $\bbT^k = \bbR^k / \bbZ^k$ is not locally convex. But one can
still talk about global convexity after taking the universal covering of
$P/\Gamma$.

Let us formulate some  global convexity theorems based on the above
discussions. First consider the case when ${P/\Gamma}$ is
simply-connected. Then ${P/\Gamma}$ can be immersed into $\bbR^k$ (where
$n$ is the dimension of $P/\Gamma$) by an integral affine map: a local
integral affine system of coordinate on ${P/\Gamma}$ can be extended
globally because there is no monodromy. Denote by $\fj: P/\Gamma
\rightarrow \bbR^k$ such an immersion (it is unique up to integral affine
automorphisms of $\bbR^k$). Then the map
\begin{equation}
\fj \circ \widehat{\mu}: \widehat{M} \rightarrow \bbR^k \ ,
\end{equation}
where $\mu$ is the momentum map of a (quasi-)Hamiltonian $\Gamma
\rightrightarrows P$ manifold $(M,\sigma)$, is a locally injective
integral affine map. When $M$ is compact then $\widehat{M}$ is also
compact, and Lemma \ref{lem:convex_affine_map1} applies. When $M$ is
noncompact, things are more subtle: the composed map $\fj \circ
\widehat{\mu}$ may be non-proper even when $\mu$ is proper, so we need an
additional condition besides the properness of $\mu$. For example, that
$\fj$ is an embedding and its image $\fj(P/\Gamma)$ closed in $\bbR^k$
(then $\fj \circ \widehat{\mu}$ will be proper provided that $\mu$ is
proper, and we can apply Lemma \ref{lem:convex_affine_map2}). Or that
$\fj$ is an embedding and $\fj(P/\Gamma)$ is convex (not necessarily
closed) in $\bbR^k$: then we can write $\fj(P/\Gamma) =
\bigcup_{k=n}^\infty K_n$ where $K_n \subset K_{n+1}$ and each $K_n$ is a
compact convex subset of $\bbR^k$, and we can use arguments absolutely
similar to the above proof of Lemma \ref{lem:convex_affine_map2} to show
that the image of $\fj \circ \widehat{\mu}$ is convex if $\mu$ is proper.
Thus we have:

\begin{thm}
\label{thm:global_convex1} Let $(M,\sigma)$ be a connected
quasi-Hamiltonian manifold of a proper quasi-symplectic groupoid
$(\Gamma \rightrightarrows P,\omega+\Omega)$ with coad-connected
isotropy groups, with a proper momentum map $\mu$. Assume that the
orbit space $P/\Gamma$ of $\Gamma$ is simply-connected, and denote
by $\fj: P/\Gamma \rightarrow \bbR^k$ an integral affine immersion
from $P/\Gamma$ to $\bbR^k$. Assume that at least one of the
following
additional conditions is satisfied: \\
1) $M$ is compact. \\
2) $\fj$ is an embedding and $\fj(P/\Gamma)$ is closed in $\bbR^k$. \\
3) $\fj$ is an embedding and $\fj(P/\Gamma)$ is convex in $\bbR^k$. \\
Then the transverse momentum map $\widehat{\mu}$ and the composed map
$\fj\circ \widehat{\mu}$ are injective, and the image $\fj \circ
\widehat{\mu}(\widehat{M}) = \fj (\mu(M)/\Gamma)$ is a convex subset in
$\bbR^k$ with locally polyhedral boundary. (We don't count boundary points
which lie in the closure of $\fj (\mu(M)/\Gamma)$ but not in $\fj
(\mu(M)/\Gamma)$). In particular, $\widehat{M}$ with its integral affine
structure is isomorphic to a convex subset of $\bbR^k$ with locally
polyhedral boundary.
\end{thm}

\Remark. In the above theorem, one may say that $(M,\sigma)$ itself is a
convex (quasi-)Hamiltonian manifold, in analogy with Knop
\cite{Knop-Convexity2002}. The space $\widehat{M}$ with its affine
structure is not only intrinsically locally convex, but it is also
intrinsically globally convex, and the convexity of $\fj (\mu(M)/\Gamma)$
is just a manifestation of this intrinsic convexity.

Consider now the case when $P/\Gamma$ is cannot be affinely immersed into
$\bbR^k$ (in particular it is not simply-connected). Then its universal
covering is a simply-connected affine manifold and hence can be immersed
into $\bbR^k$. In fact, there is a minimal connected covering of
$P/\Gamma$ which can be immersed into $\bbR^k$: it corresponds to the
subgroup of the fundamental group of $P/\Gamma$ which is precisely the
isotropy group of the monodromy representation of the flat affine
structure of $P/\Gamma$. We will denote this minimal ``flattening''
covering by $\widetilde{P/\Gamma}$, and again by $\fj:
\widetilde{P/\Gamma} \rightarrow \bbR^k$ the corresponding integral affine
immersion. The locally injective affine map $\widehat{\mu}: \widehat{M}
\rightarrow P/\Gamma$ can be lifted to a locally injective affine map
\begin{equation}
\widetilde{\widehat{\mu}}: \widetilde{\widehat{M}} \rightarrow
\widetilde{P/\Gamma} \ ,
\end{equation}
where $\widetilde{\widehat{M}}$ denotes the (connected) covering of
$\widehat{M}$ corresponding to the subgroup of the fundamental group of
$\widehat{M}$ consisting of those elements which are mapped by
$\widehat{\mu}$ into the isotropy group of the monodromy representation of
the affine structure of $P/\Gamma$. It follows from the definition that
$\widetilde{\widehat{\mu}}$ has the following property: if $x$ and $y$ are
two points in $\widetilde{\widehat{M}}$, $x \neq y$, which project to a
same point on $\widehat{M}$, then $\widetilde{\widehat{\mu}}(x) \neq
\widetilde{\widehat{\mu}}(y)$. This property implies that
$\widetilde{\widehat{\mu}}$ is a proper map, provided that $\mu$ (or
$\widehat{\mu}$) is proper. Note that $\widetilde{\widehat{M}}$ is not
compact in general, even when $M$ is compact. We have the following analog
of Theorem \ref{thm:global_convex1}:

\begin{thm}
\label{thm:global_convex2} Let $(M,\sigma)$ be a connected
quasi-Hamiltonian manifold of a proper quasi-symplectic groupoid
$(\Gamma \rightrightarrows P,\omega+\Omega)$ with coad-connected
isotropy groups, with a proper momentum map $\mu$. With the above
notations, assume that at least one of the following conditions is
satisfied: \\
1) $\widetilde{\widehat{M}}$ is compact. \\
2) $\fj: \widetilde{P/\Gamma} \rightarrow \bbR^k$ is injective and its
image is convex. \\
3) $\fj: \widetilde{P/\Gamma} \rightarrow \bbR^k$ is injective and its
image is closed. \\
Then $\fj \circ \widetilde{\widehat{\mu}}: \widetilde{\widehat{M}}
\rightarrow \bbR^k$ is an injective integral affine map and its image is a
convex subset of $\bbR^k$ with locally polyhedral boundary. (We don't
count boundary points which lie in the closure of the image but not in the
image itself).
\end{thm}

\Remark. In this paper we didn't touch the ``real'' case, i.e. the case
with an anti-symplectic involution studied by Duistermaat
\cite{Duistermaat-Convexity1983} and other people. We would conjecture
that the main results of this paper have analogs in the case with an
(anti-symplectic) involution.

\vspace{0.3cm} {\bf Acknowledgements}. The paper arises from
numerous discussions with Alan Weinstein and aims to study some of
his questions and carry out some of his ideas on proper groupoids
and convexity. The main results of the paper were first announced at
\emph{AlanDay}, a colloquium in honor of Alan Weinstein held at
Utrecht University in 03/2003. The paper is dedicated to him. I'm
also thankful to Ping Xu for telling me about quasi-symplectic
groupoids and their connections with different momentum map
theories, and to Tudor Ratiu, Eugene Lerman, and Reyer Sjamaar for
interesting discussions on convexity and connectedness properties of
momentum maps. Thanks are due to Philippe Monnier for finding out
some errors in a previous version of the paper. I would also like to
thank the referees for numerous suggestions which helped me improve
the paper.








\providecommand{\bysame}{\leavevmode\hbox to3em{\hrulefill}\thinspace}

\end{document}